\documentclass[12pt,reqno,sumlimits]{amsart}

 \usepackage{amssymb,amscd,amsmath,epsfig}
\usepackage{ 
  array,
   booktabs,
   dcolumn,
   rotating,
   shortvrb,
   tabularx,
   units,
   url,
 }
\newtheorem{thm}{Theorem}[section]
\newtheorem{cor}[thm]{Corollary}
\newtheorem{lem}[thm]{Lemma}
\newtheorem{prop}[thm]{Proposition}
\theoremstyle{definition}
\newtheorem{defn}{Definition}[section]
\newtheorem{rem}{Remark}[section]

\newtheorem{exm}[thm]{Example}

\newcommand{\R}{{\mathbb R}}

\newcommand{\Q}{{\mathbb Q}}
\newcommand{\C}{{\mathbb C}}

\newcommand{\Z}{{\mathbb Z}}

\newcommand{\g}{{\mathfrak  g}}

\newcommand{\calB}{{\mathcal B}}

\newcommand{\calE}{{\mathcal E}}

\newcommand{\calG}{{\mathcal G}}
\newcommand{\calH}{{\mathcal H}}

\newcommand{\calR}{{\mathcal R}}

\renewcommand{\to}{\longrightarrow}

\newcommand{\ev}{\operatorname{ev}}

\newcommand{\Tr}{\operatorname{Tr}}
\newcommand{\tr}{\operatorname{tr}}

\newcommand{\ad}{{\operatorname{ad\,}}}

\newcommand{\weight}{(1+\Delta)^s}
\newcommand{\weightinv}{(1+\Delta)^{-s}}

\newcommand{\con}[2]{\nabla^{#1}_{#2}}

\newsavebox{\savepar}

\newcommand{\ip}[1]{\langle #1 \rangle}

\numberwithin{equation}{section}

\newcounter{labelflag} 
\newcommand{\labelon}{\setcounter{labelflag}{1}}
\newcommand{\Label}[1]{
                       \ifnum\thelabelflag=1
                          \ifmmode
                             \makebox[0in][l]{\qquad\fbox{\rm#1}}
                          \else
                             \marginpar{\vspace{0.7\baselineskip}
                                        \hspace{-1.1\textwidth}
                                        \fbox{\rm#1}}
                          \fi
                       \fi
                       \label{#1}
                      }
\labelon

  \newcommand{\BbR}{{\mathbb R}}

 \newcommand{\BbC}{{\mathbb C}}
 \newcommand{\BbZ}{{\mathbb Z}}
 \newcommand{\pdo}{\Psi{\rm DO}}
 \newcommand{\calg}{{\mathfrak g}}
 
 \newcommand{\e}{\varepsilon}

 \newcommand{\eee}{e^{i(\theta-\theta')\cdot\xi}}
 
 \newcommand{\dg}{\dot\gamma}
 \newcommand{\ch}[3]{\Gamma_{{#1}{#2}}^{{#3}}}
 \newcommand{\cc}[3]{c_{{#1}{#2}}^{{#3}}}

 \newcommand{\ints}{\int_{S^1}}
 \newcommand{\intss}{\int_{S^*S^1}}
 \newcommand{\dtau}{\frac{\partial}{\partial\tau}\left|_{_{_{_{_{_{\tau =
 	      0}}}}}} }
 \newcommand{\ptau}[1]{\partial_\tau^{#1}}
 \newcommand{\xii}[1]{ (\xi^2)^{#1} }
 \newcommand{\eff}[1]{e^{i{#1}\theta}}

 \newcommand{\dir}{\partial\kern-.570em /}
 \newcommand{\dire}{\partial\kern-.570em /{}^{\rm eq}}
 \newcommand{\GG}{\pdo_0^*}

 \newcommand{\pa}{\partial}
 \newcommand{\cch}[2]{\Gamma_{{#1}}^{{#2}}}

 \newcommand{\xly}{(X\leftrightarrow Y)}

 \newcommand{\chw}[3]{\Gamma_{{#1}{#2}}^{{#3}} }
 \newcommand{\wgti}{(1+\Delta)^{-1}}
 \newcommand{\wgt}{1+\Delta}
 \newcommand{\ndg}{\nabla_{\dot\gamma}}
 \newcommand{\wgtsi}{(1+\Delta)^{-s}}
 \newcommand{\wgts}{(1+\Delta)^s}
 \newcommand{\ips}[2]{\langle {#1},{#2}\rangle_{s}}
 \newcommand{\ipo}[2]{\langle {#1},{#2}\rangle_{0}}

 \newcommand{\diffm}{{\rm Diff}(M) }
 \newcommand{\diff}{{\rm Diff} }
 \newcommand{\maps}{{\rm Maps} }
 \newcommand{\kk}{2k-1 }

\begin{document}

\newpage\null\vskip-4em
  \noindent\scriptsize{}
\vskip 0.5 in

\normalsize

\title{Secondary Characteristic Classes on Loop Spaces}
\author[Y. Maeda]{Yoshiaki Maeda}
\address{Department of Mathematics\\
Keio University}
\email{maeda@math.keio.ac.jp}
\author[S. Rosenberg]{Steven Rosenberg}
\address{Department of Mathematics and Statistics\\
  Boston University}
\email{sr@math.bu.edu}
\author[F. Torres-Ardila]{Fabi\'an Torres-Ardila}
\address{COSMIC\\
  University of Massachusetts Boston}
\email{fatorres@umb.edu}

\begin{abstract}   A Riemannian metric on a manifold $M$ induces a family of
  Riemannian metrics on the loop space $LM$ depending on a Sobolev space
  parameter $s$.   The connection and curvature forms of these metrics
  take values in
  pseudodifferential operators ($\Psi$DOs).
 We develop a theory of Wodzicki-Chern-Simons classes
  $CS^W_{K}\in H^{2k-1}(LM^{2k-1})$, for $K = (k_1,...,k_\ell)$ a partition of $2k-1$,
 using the the $s=0, 1$ connections and the Wodzicki residue on $\Psi$DOs. 
The new invariant $CS_5^W$ distinguishes the smooth homotopy type of
 certain actions on 
 $S^2\times S^3$, and allows us to show that $\pi_1(\diff(S^2\times S^3))$
 is infinite. \end{abstract}
\maketitle
\centerline{Dedicated to the memory of Prof. Shoshichi Kobayashi}

\bigskip\bigskip
\centerline{\sc Contents}
\medskip

\noindent 1. Introduction
\medskip

\noindent {\bf Part I. The Levi-Civita Connection on the Loop Space $LM$}
\medskip

\noindent 2. The Levi-Civita Connection for  Sobolev Parameter $s\geq 0$

\noindent 3. Local Symbol Calculations

\noindent 4. The Loop Group Case
\medskip

\noindent {\bf Part II. Characteristic Classes on $LM$}
\medskip

\noindent 5. Chern-Simons Classes on Loop Spaces



\noindent 6. An Application of Wodzicki-Chern-Simons Classes to Circle Actions

\medskip
\noindent References

\bigskip
\section{{\bf Introduction}}

The loop space $LM$ of a manifold $M$ appears frequently in mathematics and
 mathematical physics. In this paper, using an infinite dimensional version of
 Chern-Simons theory associated to the Wodzicki residue for
 pseudodifferential operators ($\pdo$s), we develop a  computable theory of secondary
 characteristic classes  on the tangent
 bundle to loop spaces.   We apply these secondary classes to distinguish circle actions
 on $S^2\times S^3$, and we prove that $\pi_1(\diff(S^2\times S^3))$ is infinite.  To our knowledge, these applications are the first examples of nonzero Wodzicki-type characteristic classes. 

Since Chern-Weil and Chern-Simons theory are geometric, it is necessary
to understand connections and curvature on loop spaces.  A  Riemannian metric
$g$ 
on $M$ induces a family of metrics $g^s$ on $LM$ parametrized by a Sobolev
space parameter $s \geq 0$, where $s=0$ gives the usual $L^2$ metric, and the
smooth case is a kind of limit as $s\to\infty.$  Thus we think of
$s$ as a regularizing parameter, and pay attention to  the parts of the theory which
are independent of $s$. 

In Part I, we compute the connection and curvature  for the Levi-Civita
connection for $g^s$ for $s>\frac{1}{2}$.  The 
closed form expressions obtained for the Levi-Civita connection for general $LM$ 
extend Freed's results for loop groups \cite{Freed}.  The connection and 
curvature forms take values in
zeroth order $\pdo$s acting on a trivial bundle over $S^1$. 
For Wodzicki-Chern-Simons classes, we only need
the principal and subprincipal symbols for these
forms, which we calculate.

In Part II, we develop a theory of Chern-Simons classes on loop spaces.
The structure group for the Levi-Civita connection for
$(LM, g^s)$ is the set of invertible zeroth order $\pdo$s, so we need
 invariant polynomials on the corresponding Lie algebra.  The naive choice is
the standard polynomials $\Tr(\Omega^k)$ of the curvature $\Omega = \Omega^s$,
where Tr is the operator trace.  However, $\Omega^k$ is zeroth order
and hence not trace class, and in any case the operator trace
is impossible to compute in general.  Instead, as in \cite{P-R2} we use the 
Wodzicki residue, the only trace on the full
algebra of $\pdo$s.  
Following Chern-Simons
 \cite{C-S} as much as possible, we build a theory of
Wodzicki-Chern-Simons (WCS) classes, which gives classes in $H^{2k-1}(LM^{2k-1})$ associated to partitions of $k$.

There are two main differences from the finite
dimensional theory. The absence of a Narasimhan-Ramanan universal connection
theorem means that we do not have a theory of differential characters
\cite{Ch-Si}.  However, since we have a family of connections on $LM$, we can define real valued, not just $\R/\Z$-valued, WCS classes. 

 In contrast to the operator trace, the Wodzicki residue is locally
 computable, so we can write explicit expressions for the WCS classes.
In particular, we can see how the WCS classes depend on the Sobolev parameter $s$, and 
hence define a  ``regularized" or $s$-independent WCS classes.  
The local expression also yields some vanishing results for WCS classes.  More importantly,
we produce a nonvanishing 
WCS class on $L(S^2\times S^3).$  This leads to the topological results described in the first paragraph.

For related results on characteristic classes on infinite rank bundles with a group of $\pdo$s as structure group, see 
 \cite{lrst, P-R2}.

\medskip
The paper is organized as follows.  Part I treates the family of metrics $g^s$ on $LM$
associated to $(M,g)$.  \S2 discusses connections associated to $g^s.$  After some preliminary material,
we compute the Levi-Civita connection for $s=0$ (Lemma \ref{lem:l2lc}), $s=1$
(Theorem \ref{old1.6}),  $s\in \Z^+$ (Theorem \ref{thm:sinz}), and general 
$s>\frac{1}{2}$ (Theorem \ref{thm25}).   These connections allow us to track how the geometry
of $LM$ depends on $s$.  

Both the Levi-Civita and $H^s$ connections have connection and curvature forms taking 
values in $\pdo$s of order zero.  In \S3, we compute the symbols of these forms needed in Part II.  In \S4, we show that our results extend Freed's on loop groups \cite{Freed}.

Part II covers Wodzicki-Chern-Simons classes.  In \S5, we review the finite dimensional
construction of Chern and Chern-Simons classes, and use the Wodzicki residue to define Wodzicki-Chern (WC) and WCS classes (Definition \ref{def:WCS}).  We prove the necessary vanishing
of the WC classes for mapping spaces (and in particular for $LM$) in Proposition
\ref{prop:maps}.  In Theorem \ref{thm:5.5}, we give the explicit local expression for the
relative WCS class $CS_{2k-1}^W(g)\in H^{2k-1}(LM^{2k-1})$ associated to the trivial
partition of $k$.  We then define the regularized or $s$-independent WCS class.  
In Theorem \ref{WCSvan}, we give a vanishing result
for WCS classes.  

In particular, the WCS class which is the analogue of the classical dimension three Chern-Simons class vanishes on loop spaces of 
 $3$-manifolds, so we look for nontrivial examples on $5$-manifolds.
 In \S\ref{dimfive}, we use a Sasaki-Einstein
metric constructed in \cite{gdsw} to produce a nonzero WCS class  $CS_5^W\in
H^5(L(S^2\times
S^3)).$  We prove $CS_5^W\neq 0$ by an exact computer calculation showing
 $\int_{[a^L]} CS_5^W \neq 0$, where
$[a^L]\in H_5(LM)$ is a cycle associated to a simple
circle action on $S^2\times S^3.$  From this
nonvanishing, we conclude both that  the circle action is not 
smoothly homotopic to the trivial action and that $\pi_1(\diff(S^2\times S^3))$ is infinite.
We expect other similar results in the future.

Our many discussions with Sylvie Paycha are gratefully
acknowledged. We also thank Kaoru Ono and Dan Freed for pointing out errors in  previous versions of 
the paper.

\bigskip

\large
\noindent {{\bf Part I. The Levi-Civita Connection on the Loop Space $LM$}}
\normalsize

\bigskip

In this part, we compute the Levi-Civita connection on $LM$
associated to a Riemannian metric on $M$ and a Sobolev parameter $s=0$ 
or $s>\frac{1}{2}.$  The standard $L^2$ metric on $LM$ is the case $s=0$, and otherwise we avoid technical issues by assuming that $s$ is greater than the critical exponent $\frac{1}{2}$ for analysis on bundles over $S^1.$  
 In \S\ref{LCconnection}, the
main results are  Lemma \ref{lem:l2lc}, Theorem \ref{old1.6}, 
Theorem \ref{thm:sinz}, and Theorem \ref{thm25},
which compute the Levi-Civita connection for $s =0$, $s=1$, 
$s\in \Z^+$,  and general $s >\frac{1}{2},$ respectively.

In
\S3, we compute the relevant symbols of the connection one-forms and the
curvature two-forms.  In \S4, we compare our results with work of Freed
\cite{Freed} on loop groups.

\section{{\bf The Levi-Civita Connection for Sobolev Parameter $s\geq 0$}}
\label{LCconnection}

This section covers background material and computes the Levi-Civita
connection on $LM$ for Sobolev parameter $s=0$ and $s>\frac{1}{2}$.  
In \S2.1, we review
material on $LM$, and in \S2.2 we review pseudodifferential operators and the
Wodzicki residue.  In \S2.3, we give the crucial computations of the Levi-Civita connections
for $s=0,1$.
This computation is extended to $s\in \Z^+$ in \S2.4, and to general $s>\frac{1}{2}$ in
\S2.5.    In \S2.6, we discuss how the geometry of $LM$ forces an extension of
the structure group of $LM$ from a gauge group to a group of bounded
invertible $\pdo$s.
     
\subsection{{\bf Preliminaries on $LM$}}

${}$
\medskip

Let $(M, \langle\ ,\ \rangle)$ 
be a closed, connected,  oriented Riemannian $n$-manifold with  loop space $LM
= C^\infty(S^1,M)$ of smooth loops. 
$LM$ is a smooth infinite dimensional Fr\'echet manifold, but it is
 technically simpler 
to work 
with the smooth Hilbert manifold $H^{s'}(S^1,M)$ of loops in some Sobolev class $s' \gg 0,$
as we now recall. For $\gamma\in LM$, the formal
tangent space $T_\gamma LM$ is 
$\Gamma(\gamma^*TM)$, the space
 of smooth sections of the pullback bundle $\gamma^*TM\to
S^1$.   The actual tangent space of $H^{s'}(S^1, M)$ at $\gamma$ is 
$H^{s'-1}(\gamma^*TM),$     the sections of $\gamma^*TM$ of Sobolev class $s'-1.$
We will fix $s'$ and use $LM, T_\gamma LM$ for $H^{s'}(S^1, M), H^{s'-1}(\gamma^*TM)$, respectively.

For each $s>1/2,$ we can complete $\Gamma(\gamma^*TM\otimes \C)$
 with respect to the Sobolev inner product
 \begin{equation}\label{eq:Sob1}
\langle X,Y\rangle_{s}=\frac{1}{2\pi}\int_0^{2\pi} \langle(1+\Delta)^{s}
X(\alpha),Y(\alpha)
\rangle_{\gamma (\alpha)}d\alpha,\  X,Y\in \Gamma(\gamma^*TM).
\end{equation}
Here $\Delta=D^*D$, with $D=D/d\gamma$ the covariant derivative along
$\gamma$. (We use this notation instead of the classical $D/dt$ to keep track
of $\gamma$.)
We need the complexified pullback bundle $\gamma^*TM\otimes \BbC$, denoted from now on
just as $\gamma^*TM$, in order to apply the
pseudodifferential operator  $(1+\Delta)^{s}.$
The construction of $(1+\Delta)^{s}$ is reviewed in
\S\ref{pdoreview}.  
We denote this completion by  $H^{s'}(\gamma^*TM)$.   We can consider the
$s$ metric on $TLM$ for any $s\in \R$, but we will only consider 
$s=0$ or $1/2 < s\leq s'-1.$

A small real neighborhood $U_\gamma$ 
of the zero section in $H^{s'}(\gamma^*TM)$ is a
coordinate chart near $\gamma\in LM$ 
via the pointwise exponential map
\begin{equation}\label{pointwiseexp}
\exp_\gamma:U_\gamma
\to L M, \ X \mapsto 
\left(\alpha\mapsto \exp_{\gamma(\alpha)} X(\alpha)\right).  
\end{equation}
Note that the domain of the exponential map is not contained in $T_\gamma LM.$
The differentiability of the transition functions $\exp_{\gamma_1}^{-1}\cdot
\exp_{\gamma_2}$ is proved in
\cite{E} and \cite[Appendix A]{Freed1}.
Here $\gamma_1, \gamma_2$ are close loops in the sense that
a geodesically convex neighborhood of $\gamma_1(\theta)$ contains
$\gamma_2(\theta)$ and vice versa for all $\theta.$
Since 
$\gamma^*TM$
is (noncanonically) isomorphic to the trivial bundle ${\mathcal R} =
S^1\times \BbC^n\to S^1$, 
the model space for $LM$ is the set of 
$H^{s'}$ sections of this trivial bundle.  
The $s$ metric is a weak Riemannian metric for $s<s'$ in the sense that the topology induced on $H^{s'}(S^1, M)$ by the exponential map applied to $H^{s}(\gamma^*TM)$ is weaker than the $H^{s'}$ topology.

The complexified tangent bundle
$TLM$ has
transition functions $d(\exp_{\gamma_1}^{-1} 
\circ \exp_{\gamma_2})$.  Under
the isomorphisms $\gamma_1^*TM \simeq {\mathcal R} \simeq
\gamma_2^*TM$, the transition functions lie in the gauge group
${\calG}({\mathcal R})$, so this is the structure group of $TLM.$


\subsection{{\bf Review of $\pdo$ Calculus}}\label{pdoreview}
${}$
\medskip

We recall the construction of classical
pseudodifferential operators ($\pdo$s)
 on a closed 
manifold $M$ from \cite{gilkey, See}, assuming knowledge of $\pdo$s on
 $\R^n$ (see e.g. \cite{hor, shu}).   
 

 A linear operator $P:C^\infty(M)\to C^\infty(M)$
is a $\pdo$ of  order $d$ if for
every open chart $U\subset M$ and functions $\phi,\psi\in C_c^\infty(U)$,
$\phi P\psi$ is a $\pdo$ of order $d$ on $\R^n$, where we do not
 distinguish between $U$ and its diffeomorphic image in $\R^n$.
 Let
  $\{U_i\}$ be a finite cover of $M$
with subordinate partition of unity
  $\{\phi_i\}.$  Let $\psi_i\in C^\infty_c(U_i)$ have $\psi_i \equiv 1$ on
  supp$(\phi_i)$  and set $P_i = \psi_iP\phi_i.$
Then
  $ \sum_i \phi_iP_i\psi_i$ is a $\pdo$ on $M$, and
$P$ differs from $ \sum_i \phi_iP_i\psi_i$ by a smoothing operator, denoted
$P\sim \sum_i \phi_iP_i\psi_i$. 
 In particular, this sum is independent of the
  choices up to smoothing operators.
All this carries over to $\pdo$s acting on sections of a bundle over $M$.

An example is the $\pdo$ $(1+\Delta-\lambda)^{-1}$ for $\Delta$ a positive
order nonnegative elliptic $\pdo$ and $\lambda$
outside the spectrum of $1+\Delta.$  In each $U_i$, we construct a parametrix
$P_i$ for $A_i = \psi_i
(1+\Delta-\lambda)\phi_i$ by formally inverting $\sigma(A_i)$ and then
constructing a $\pdo$ with the inverted symbol.  By \cite[App.~A]{A-B2},
$ B = \sum_i \phi_iP_i\psi_i$ is a parametrix for $(1+\Delta-\lambda)^{-1}$.
  Since
  $B \sim (1+\Delta-\lambda)^{-1}$,
  $(1+\Delta-\lambda)^{-1}$ is itself a $\pdo$.
For $x\in U_i$, by definition
$$\sigma((1+\Delta-\lambda)^{-1})(x,\xi) = \sigma(P)(x,\xi) = \sigma(\phi
P\phi)(x,\xi),$$
where $\phi$ is a bump function with $\phi(x) = 1$ \cite[p.~29]{gilkey};  
the symbol depends on the choice of $(U_i, \phi_i).$  

The operator $\weight$ for Re$(s) <0$,
which exists as a bounded
operator on $L^2(M)$ by the functional calculus, is also a $\pdo$.
To see this,
 we construct the putative symbol $\sigma_i$
of $ \psi_i\weight\phi_i$ in each $U_i$
by a contour integral $\int_\Gamma \lambda^{s}\sigma[(1+\Delta-\lambda)^{-1}]
d\lambda$
around the spectrum of $1+\Delta$.
We then
construct a $\pdo$ $Q_i$ on $U_i$ with $\sigma(Q_i) = \sigma_i$,
and set $Q  = \sum_i\phi_i Q_i\psi_i.$
By arguments in \cite{See},
$\weight \sim Q$,
so $\weight$ is a $\pdo$.


 Recall that the {\it Wodzicki residue} of a $\pdo$ $P$ on sections
  of a bundle $E\to
M^n$ is
\begin{equation}\label{wresdef} {\rm res}^{\rm w}(P) = 
\int_{S^*M} \tr\ \sigma_{-n}(P)(x,\xi) d\xi dx,\end{equation}
where $S^*M$ is the unit cosphere bundle for some metric.  The Wodzicki
residue is independent of choice of local coordinates, and up to scaling is 
the unique trace on the algebra of $\pdo$s if dim$(M)>1$ (see
e.g.~\cite{fgl} in general and \cite{ponge} for the case $M=S^1.$).

The Wodzicki residue will be used in
Part II to define characteristic classes on $LM$.  In our particular case, the operator $P$ 
will be an $\pdo$ of order $-1$ acting on sections of a bundle over $S^1$
 (see (\ref{cswint})), so
$\sigma_{-1}(P)$ is globally defined. Of course, $\int_{S^*S^1}
\tr\sigma_{-1}(P) d\xi d\theta = 2\int_{S^1}\tr\sigma_{-1}(P)d\theta$.  It is
easy to check that this integral,
which strictly speaking involves a choice of cover of $S^1$ and a partition of unity,
equals the usual $2\int_0^{2\pi} \tr\sigma_{-1}(P) d\theta.$



\subsection{The Levi-Civita Connection for $s=0, 1$}
${}$
\medskip

The smooth Riemannian manifold $LM = H^{s'}(S^1,M)$ has tangent bundle $TLM$ with 
$T_\gamma LM = H^{s'-1}(\gamma^*TM).$  For the $s'-1$ metric on $TLM$ (i.e., 
$s = s'-1$ in (\ref{eq:Sob1})), 
the 
Levi-Civita connection exists and is determined by the six term formula
\begin{eqnarray}\label{5one}
2\ip{\con{s}{X}Y,Z}_{s} &=& X\ip{Y,Z}_{s}+Y\ip{X,Z}_{s}-Z\ip{X,Y}_{s}\\
&&\qquad +\ip{[X,Y],Z}_{s}+\ip{[Z,X],Y}_{s}-\ip{[Y,Z],X}_s\nonumber
\end{eqnarray}
\cite[Ch. VIII]{lang}.  The point is that each term on the RHS of (\ref{5one}) 
is
a {\it continuous} linear functional $T_i:H^{s=s'-1}(\gamma^*TM) \to \BbC$ in $Z$.  Thus 
$T_i(Z) = \ip{T_i'(X,Y),Z}_s$ for a unique $T'(X,Y)\in H^{s'-1}(\gamma^*TM)$, and $\con{s}{Y}X
= \frac{1}{2}\sum_i T'_i.$  

In general, the Sobolev parameter $s$ in (\ref{eq:Sob1}) differs from the parameter $s'$ defining the loop space.  We discuss how this affects the existence of a Levi-Civita connection. 

\begin{rem}\label{lcrem}  For general $s >\frac{1}{2}$, the Levi-Civita connection for the $H^s$ 
metric is guaranteed to exist on the bundle $H^s(\gamma^*TM)$, as above.  However, it is inconvenient to have the bundle depend on the Sobolev parameter, for several reasons:  
(i) $H^s(\gamma^*TM)$ is strictly speaking not the tangent bundle of $LM$, (ii) for the
$L^2$ ($s=0$) metric, the Levi-Civita connection should be given by the Levi-Civita connection on $M$ applied pointwise along the loop (see Lemma \ref{lem:l2lc}), and on $L^2(\gamma^*TM)$  this would have to be interpreted in the distributional sense; (iii) to compute Chern-Simons classes on
$LM$ in Part II, we need to compute with a pair of connections corresponding to $s=0, s=1$ on the
same bundle.  These problems are not fatal: (i) and (ii) are essentially  aesthetic issues,
and for (iii), the connection one-forms will take values in zeroth order $\pdo$s, which are bounded operators on any 
$H^{s'-1}(\gamma^*TM)$, so $s' \gg 0$ can be fixed.  

Thus it is more convenient
to fix $s'$ and consider the family of $H^s$ metrics on $TLM$ for 
$\frac{1}{2} < s < s'-1$. 
  However, the existence of the Levi-Civita connection for the $H^s$ metric is trickier.
  For a sequence $Z\in H^{s'-1} = H^{s'-1}(\gamma^*TM)$ with $Z\to 0$
  in $H^{s'-1}$ or in 
  $H^s$, the RHS of (\ref{5one}) goes to $0$ for fixed $X, Y\in H^s.$  Since
  $H^{s'-1}$ is dense in $H^{s}$, the RHS of (\ref{5one}) extends to a continuous linear functional on $H^s$.  Thus the RHS of (\ref{5one}) is given by
  $\langle L(X,Y), Z\rangle_s$ for some $L(X,Y)\in H^s.$  We set $\nabla^{s}_YX = 
  \frac{1}{2}L(X,Y)$.  Note that even if we naturally demand that
  $X, Y\in H^{s'-1}$, we only get $\nabla^s_YX\in H^s\supset H^{s'-1}$ without additional work.  Part of the content of Theorem \ref{thm25} is that the Levi-Civita connection exists in the {\it strong sense}:  given a tangent vector $X\in H^{s'-1}(\gamma^*TM)$ and a smooth vector field
$Y(\eta)\in H^{s'-1}(\eta^*TM)$ for all $\eta$,
   $\nabla^s_XY(\gamma)\in H^{s'-1}(\gamma^*TM).$  See Remark 2.6.
  

\end{rem}

We need to discuss local coordinates on $LM$.
For motivation, recall that
\begin{equation}\label{lie}[X,Y]^a = X(Y^a)\partial_a - Y(X^a)\partial_a
\equiv \delta_X(Y) -\delta_Y(X)
\end{equation}
in local coordinates on a finite dimensional manifold.  Note that
$X^i\partial_iY^a = X(Y^a) =
(\delta_XY)^a$ in this notation.

Let $Y$ be a vector field on $LM$, and let $X$ be a tangent vector at
$\gamma\in LM.$  The local variation $\delta_XY$ of $Y$ in the direction of $X$ at $\gamma$ is 
defined as usual: let $\gamma(\e,\theta)$ be a family of loops in $M$
with $\gamma(0,\theta) = \gamma(\theta), \frac{d}{d\e}|_{_{\e=0}}
\gamma(\e,\theta) = X(\theta).$  Fix $\theta$, and let $(x^a)$ be
coordinates near $\gamma(\theta)$.  We call these coordinates 
{\it manifold coordinates.} Then
$$\delta_XY^a(\gamma)(\theta) \stackrel{{\rm def}}{=}
\frac{d}{d\e}\biggl|_{_{_{\e =0}}} Y^a(\gamma(\e,\theta)).$$
Note that $\delta_XY^a = (\delta_XY)^a$ by definition.

\begin{rem} Having $(x^a)$ defined only near a fixed $\theta$ is inconvenient.
We can find coordinates that work for all $\theta$ as follows. For
  fixed $\gamma$, there is an $\e$ such that for all $\theta$,
  $\exp_{\gamma(\theta)} X$ is inside the cut locus of $\gamma(\theta)$ if
  $X\in T_{\gamma(\theta)}M$ has $|X|<\e.$  Fix such an $\e.$ Call 
  $X\in H^{s'-1}(\gamma^*TM)$ {\it
  short} if $|X(\theta)|<\e$ for all $\theta.$  Then
$$U_\gamma = \{\theta \mapsto \exp_{\gamma(\theta)}X(\theta) | X\ {\rm is\
    short}\}\subset LM$$
is a coordinate neighborhood of $\gamma$ parametrized by $\{ X: X\  {\rm is\ 
    short}\}.$  

 We know
$H^{s'-1}(\gamma^*TM) \simeq H^{s'-1}(S^1\times \R^n)$ noncanonically, so
$U_\gamma$ is parametized by short sections of $H^{s'-1}(S^1\times \R^n)$ for
a different $\e.$  In particular, we have a smooth diffeomorphism $\beta$ from
$U_\gamma$ to short sections of $H^{s'-1}(S^1\times \R^n)$.

Put coordinates $(x^a)$ on $\R^n$, which we identify canonically with the fiber $\R^n_\theta$
over $\theta$ in $S^1\times \R^n$. For $\eta\in
U_\gamma$, we have $\beta(\eta) = (\beta(\eta)^1(\theta),...,\beta(\eta)^n(\theta)).$
As with finite dimensional coordinate systems, we will drop $\beta$ and just
write
$\eta = (\eta(\theta)^a).$ These coordinates work for all
$\eta$ near $\gamma$ and for all $\theta.$  The definition of $\delta_XY$ above carries over to exponential coordinates.

We will call these coordinates {\it exponential coordinates}.
\end{rem}

(\ref{lie}) continues to hold
 for vector fields on $LM$, in either 
 manifold or exponential coordinates.
   To see this, one checks that the coordinate-free proof that $L_XY(f) =
 [X,Y](f)$ for $f\in C^\infty(M)$ (e.g.~\cite[p.~70]{warner}) carries over to
 functions on $LM$.  In brief, the usual proof involves a map $H(s,t)$ of a
 neighborhood of the origin in $\R^2$ into $M$, where $s,t$ are parameters for
 the flows of $X, Y,$ resp.  For $LM$, we have a map $H(s,t,\theta)$, where
 $\theta$ is the loop parameter.  
The  usual proof uses
 only $s, t$ differentiations,
 so $\theta$ is unaffected.    The point is that the $Y^i$ are local functions
 on the $(s,t,\theta)$  parameter space, whereas
the $Y^i$ are not
  local functions on $M$ at points where loops cross or self-intersect.

We first compute the $L^2$ ($s=0$) Levi-Civita connection invariantly and in 
manifold coordinates.

\begin{lem} \label{lem:l2lc} Let $\nabla^{LC}$ be the Levi-Civita connection on $M$.  
 Let $\ev_\theta:LM\to M$ be $\ev_\theta(\gamma) = \gamma(\theta).$ 
 Then $D_XY(\gamma)(\theta) \stackrel{\rm def}{=} 
 (\ev_\theta^*\nabla^{LC})_XY(\gamma)(\theta)$ is the $L^2$ Levi-Civita connection on $LM$.  In manifold coordinates,
 \begin{equation}\label{l2lc} (D_XY)^a(\gamma)(\theta) = \delta_XY^a(\gamma)(\theta) +
  \cch{bc}{a}(\gamma(\theta))X^b(\gamma)(\theta) Y^c(\gamma)(\theta).
  \end{equation}
\end{lem}
\medskip

As in Remark \ref{lcrem}, we may assume that
$X, Y\in H^{s'-1}(\gamma^*TM)$ with $s' \gg 0$, so (\ref{l2lc}) makes sense.

\begin{proof} $\ev_\theta^*\nabla^{LC}$ is a connection on
$\ev_\theta^*TM\to LM$. We have
$\ev_{\theta,*}(X) = X(\theta)$.  If $U$ is a coordinate
neighborhood on $M$ near some $\gamma(\theta)$, then on $\ev_\theta^{-1}(U)$, 
\begin{eqnarray*}(\ev_\theta^*\nabla^{LC})_XY^a(\gamma)(\theta) &=& (\delta_{X}Y)^a(\gamma)(\theta) +
((\ev_\theta^*\omega^{LC}_{X})Y)^a (\theta)\\
&=& (\delta_{X}Y)^a(\gamma)(\theta) + 
  \chw{b}{c}{a}(\gamma(\theta))X^b(\gamma)(\theta) Y^c(\gamma)(\theta)
\end{eqnarray*}
Since $\ev_\theta^*\nabla^{LC}$ is a connection, for each fixed $\theta$, $\gamma$ and $X\in
 T_\gamma LM$, 
 $Y\mapsto$\\
 $ (\ev^*_\theta\nabla^{LC})_XY(\gamma)$
 has Leibniz rule with respect to
functions on $LM$.   Thus $D$ is a connection on $LM.$

$D$ is torsion free, as from the local expression
  $D_XY - D_YX = \delta_XY - \delta_YX = [X,Y].$

To show that  $D_XY$ is compatible with the $L^2$
  metric, first recall that  for a function $f$ on $LM$, $D_Xf = \delta_Xf =
  \frac{d}{d\e}|_{_{\e=0}}f(\gamma(\e,\theta))$ for $X(\theta)
   = \frac{d}{d\e}|_{_{\e=0}}\gamma(\e, \theta).$
  (Here $f$ depends only on
  $\gamma$.)  Thus (suppressing the partition of unity, which is independent of $\e$)
\begin{eqnarray*} D_X\langle Y,Z\rangle_0 &=& 
  \frac{d}{d\e}\biggl|_{_{_{\e=0}}}\int_{S^1} g_{ab}(\gamma(\e,\theta))
Y^a(\gamma(\e,\theta))Z^b(\gamma(\e,\theta))d\theta\\
&=& \int_{S^1}\partial_c
g_{ab}(\gamma(\e,\theta))
X^cY^a(\gamma(\e,\theta))Z^b(\gamma(\e,\theta))d\theta\\
&&\qquad + \int_{S^1} 
g_{ab}(\gamma(\e,\theta))
(\delta_XY)^a(\gamma(\e,\theta))Z^b(\gamma(\e,\theta))d\theta\\
&&\qquad 
+ \int_{S^1} 
g_{ab}(\gamma(\e,\theta))
Y^a(\gamma(\e,\theta))(\delta_XZ)^b(\gamma(\e,\theta))d\theta\\
&=& \int_{S^1}\Gamma{}_{c a}^{e}g_{eb} X^cY^aZ^b +
\Gamma{}_{c b}^{e}g_{ae}X^cY^aZ^b\\
&&\qquad +g_{ab}(\delta_XY)^aZ^b + g_{ab}Y^a(\delta_X Z)^bd\theta\\
&=& \langle D_XY,Z\rangle_0 + \langle Y, D_XZ\rangle_0.
\end{eqnarray*}
\end{proof}

\begin{rem} The local expression for $D_XY$ also holds in exponential coordinates. More precisely, let $(e_1(\theta),...,e_n(\theta))$
be a global frame of $\gamma^*TM$ given by the trivialization of
$\gamma^*TM.$  Then $(e_i(\theta))$ is also naturally a frame of
$T_XT_{\gamma(\theta)}M$ for all $X\in T_{\gamma(\theta)}M.$  We use
$\exp_{\gamma(\theta)}$ to pull back the metric on $M$ to a metric on
$T_{\gamma(\theta)}M$: 
$$g_{ij}(X) = (\exp^*_{\gamma(\theta)}g)(e_i, e_j) =
  g(d(\exp_{\gamma(\theta)})_X (e_i),   d(\exp_{\gamma(\theta)})_X
  (e_j))_{\exp_{\gamma(\theta)}X}.$$
Then the Christoffel symbols  
$\Gamma_{b c}^{a}(\gamma(\theta))$
are computed with respect to
  this metric.  For example, the term $\partial_\ell g_{bc}$ means $e_\ell
  g(e_a, e_b)$, etc.  The proof that $D_XY$ has the local expression (\ref{l2lc}) 
  then carries over to exponential coordinates.

\end{rem}

The  $s=1$ Levi-Civita connection on $LM$ is given as follows.

\begin{thm} \label{old1.6}
The $s=1$ Levi-Civita connection $\nabla = \nabla^1$ on $LM$ is given at the loop
$\gamma$ by
\begin{eqnarray*}  \nabla_XY &=& D_XY + \frac{1}{2}\wgti\left[
-\ndg(R(X,\dg)Y) - R(X,\dg)\ndg Y\right.\\
&&\qquad \left. -\ndg(R(Y,\dg)X) - R(Y,\dg)\ndg X\right.\\
&&\qquad \left. +R(X,\ndg Y)\dg - R(\ndg X, Y)\dg\right].
\end{eqnarray*}
\end{thm}

We prove this in a series of steps.  The assumption in the next Proposition will be dropped later.

\begin{prop} \label{old1.3}
The Levi-Civita connection for the $s=1$ metric is given by
$$\nabla_X^1Y = D_XY + \frac{1}{2}\wgti[D_X, 1+\Delta]Y +
  \frac{1}{2}\wgti[D_Y, 1+\Delta]X 
+ A_XY,$$
where  we assume that for $X, Y\in H^{s'-1}$, $A_XY$ is well-defined by
\begin{equation}\label{insert2}-\frac{1}{2}\langle [D_Z,\wgt]X,Y\rangle_0 = \langle A_XY,Z\rangle_1.
\end{equation}
\end{prop}

\begin{proof} By Lemma \ref{lem:l2lc},
\begin{eqnarray*} X\langle Y,Z\rangle_1 &=& X\langle (\wgt)Y,Z\rangle_0 =
  \langle D_X((\wgt)Y),Z\rangle_0 + \langle (\wgt)Y, D_XZ\rangle_0\\
Y\langle X,Z\rangle_1 &=& \langle D_Y((\wgt)X),Z\rangle_0 + \langle (\wgt)X,
D_YZ\rangle_0\\
-Z\langle X,Y\rangle_1 &=& -\langle D_Z((\wgt)X),Y\rangle_0 - \langle (\wgt)X,
D_ZY\rangle_0\\
\langle [X,Y],Z\rangle_1 &=& \langle(\wgt)(\delta_XY - \delta_YX), Z\rangle_0
= \langle (\wgt)(D_XY - D_YX),Z\rangle_0\\
\langle[Z,X],Y\rangle_1 &=& \langle(\wgt)(D_ZX-D_XZ),Y\rangle_0\\
-\langle[Y,Z],X\rangle_1 &=& -\langle(\wgt)(D_YZ-D_ZY),X\rangle_0.
\end{eqnarray*}
The six terms on the left hand side must sum up to $2\langle \nabla^1_XY,Z\rangle_1$ 
in the sense of Remark \ref{lcrem}.
After some cancellations, for $\nabla =\nabla^1$ we get
\begin{eqnarray*} 2\langle\nabla_XY,Z\rangle_1 &=& \langle D_X((\wgt)Y),Z\rangle_0 +
  \langle D_Y((\wgt)X),Z\rangle_0\nonumber\\
&&\qquad + \langle (\wgt)(D_XY - D_YX),Z\rangle_0 - \langle
  D_Z((\wgt)X),Y\rangle_0\nonumber\\
&&\qquad +\langle(\wgt)D_ZX),Y\rangle_0\nonumber\\
&=& \langle (\wgt)D_XY,Z\rangle_0 + \langle [D_X,\wgt] Y, Z\rangle_0\\
&&\qquad + \langle (\wgt)D_YX,Z\rangle_0 + \langle [D_Y,\wgt] X, Z\rangle_0\nonumber\\
&&\qquad + \langle (\wgt)(D_XY - D_YX),Z\rangle_0 -\langle
  [D_Z,\wgt]X,Y\rangle_0\\
&=& 2\langle D_XY,Z\rangle_1 + \langle \wgti[D_X,\wgt]Y,Z\rangle_1\\
&&\qquad + \langle \wgti[D_Y,\wgt]X,Z\rangle_1 +2
\langle A_XY,Z\rangle_1.
\end{eqnarray*}

\end{proof}

Now we compute the bracket terms in the Proposition.  We have $[D_X,\wgt] =
[D_X,\Delta]$. Also,
$$0 = \dot\gamma\langle X, Y\rangle_0 = \langle\nabla_{\dot\gamma}X,Y\rangle_0
+ \langle X,\nabla_{\dot\gamma}Y\rangle_0,$$
so 
\begin{equation}\label{one}\Delta = \nabla_{\dot\gamma}^* \nabla_{\dot\gamma}
  = -\nabla_{\dot\gamma}^2.
\end{equation}

\begin{lem} $[D_X,\nabla_{\dot\gamma}]Y = R(X,\dot\gamma)Y.$
\end{lem}

\begin{proof} 
Note that $\gamma^\nu, \dot\gamma^\nu$ are locally defined functions on 
$S^1\times LM.$
Let $\tilde\gamma:
[0,2\pi]\times (-\e,\e)\to M$ be a smooth map with $\tilde\gamma(\theta,0) =
\gamma(\theta)$, and
$\frac{d}{d\tau}|_{\tau = 0}\tilde\gamma(\theta,\tau) = Z(\theta).$
Since $(\theta,\tau)$ are coordinate functions on
$S^1\times (-\e,\e)$, we have
\begin{eqnarray}\label{badterms} Z(\dg^\nu) &=& \delta_Z(\dg^\nu) = \ptau{Z}(\dg^\nu) =
\dtau\right.\left(\frac{\partial}{\partial\theta}
(\tilde\gamma(\theta,\tau)^\nu\right)\\
&=& \frac{\partial}{\partial\theta}
\dtau\right. \tilde\gamma(\theta,\tau)^\nu = \partial_\theta Z^\nu \equiv
\dot Z^\nu.\nonumber
\end{eqnarray}

 We compute
\begin{eqnarray*} 
(D_X\nabla_{\dg} Y)^a 
&=& \delta_X(\nabla_{\dg} Y)^a +
  \chw{b}{c}{a}X^b\nabla_{\dg} Y^c\\
&=& \delta_X(\dg^j\partial_jY^a + \chw{b}{c}{a}\dg^bY^c)
+ \chw{b}{c}{a}X^b(\dg^j\partial_jY^c + \chw{e}{f}{c}\dg^e Y^f)\\
&=& \dot X^j\partial_jY^a + \dg^j\partial_j\delta_XY^a +
  \partial_m\chw{b}{c}{a}X^m\dg^bY^c
+ \chw{b}{c}{a}\dot X^bY^c + \chw{b}{c}{a}\dg^b\delta_XY^c\\
&&\qquad 
+ \chw{b}{c}{a}X^b\dg^j\partial_jY^c +
  \chw{b}{c}{a}\chw{e}{f}{c}X^b\dg^eY^f.\\
(\nabla_{\dg} D_XY)^a &=& \dg^j(\partial_j(D_XY)^a +
  \chw{b}{c}{a}\dg^b (D_XY)^c)\\
&=& \dg^j\partial_j(\delta_XY^a + \chw{b}{c}{a}X^b Y^c) 
+ \chw{b}{c}{a}\dg^b(\delta_XY^c + \chw{s}{f}{c}X^eY^f)\\
&=& \dg^j\partial_j\delta_XY^a + \dg^j\partial_j\chw{b}{c}{a}X^bY^c +
  \chw{b}{c}{a}\dot X^bY^c
+ \chw{b}{c}{a}X^b\dot Y^c + \chw{b}{c}{a}\dg^b\delta_XY^c \\
&&\qquad + \chw{b}{c}{a}\chw{e}{f}{c}\dg^b X^eY^f.
\end{eqnarray*}
Therefore
\begin{eqnarray*} (D_X\nabla_{\dg}Y - \nabla_{\dg}D_XY)^a &=& \partial_m
  \chw{b}{c}{a}X^m\dg^bY^c - \partial_j \chw{b}{c}{a}\dg^j X^bY^c
  + \chw{b}{c}{a}\chw{e}{f}{c}X^b\dg^e Y^f   \\
&&\qquad 
-\chw{b}{c}{a}\chw{e}{f}{c}\dg^b X^e Y^f \\
&=& (\partial_j \Gamma_{bc}^{a} - \partial_b \chw{j}{c}{a} +\chw{j}{e}{a}\chw{b}{c}{e}-
\chw{b}{e}{a}\chw{j}{c}{e})\dg^b X^j Y^c \\
&=& R_{jbc}^{\ \ \ a}X^j\dg^b Y^c,
\end{eqnarray*}
so 
$$D_X\nabla_{\dg}Y - \nabla_{\dg}D_XY =  R(X,\dg)Y.$$
\end{proof}

\begin{cor}\label{cor:zero}
 At the loop $\gamma$, $[D_X,\Delta]Y = -\nabla_{\dg}(R(X,\dot\gamma)Y) -
  R(X,\dg)\nabla_{\dg}Y.$  In particular, $[D_X,\Delta]$ is a zeroth order
  operator.
\end{cor}

\begin{proof}  
\begin{eqnarray*} [D_X,\Delta]Y &=& (-D_X\ndg\ndg + \ndg\ndg D_X)Y \\
&=& -(\ndg D_X\ndg Y+ R(X,\dg)\ndg Y) +\ndg\ndg D_XY\\
&=& -(\ndg\ndg D_XY + \ndg(R(X,\dg)Y) + R(X,\dg)\ndg Y) 
+\ndg\ndg D_XY\\
&=& -\ndg(R(X,\dg)Y) - R(X,\dg)\ndg Y.
\end{eqnarray*}
\end{proof}

Now we complete the proof of Theorem \ref{old1.6}, showing in the process that $A_XY$ exists. 
\medskip

\noindent {\it Proof of Theorem \ref{old1.6}.}
 By Proposition \ref{old1.3} and Corollary \ref{cor:zero}, we have
\begin{eqnarray*} \nabla_XY &=& D_XY + \frac{1}{2}\wgti[D_X,\wgt]Y +
  (X\leftrightarrow Y) + A_XY\\
&=& D_XY + \frac{1}{2}\wgti(-\ndg(R(X,\dg)Y) - R(X,\dg)\ndg Y) + 
  (X\leftrightarrow Y) + A_XY,
\end{eqnarray*}
where $(X\leftrightarrow Y)$ denotes the previous term with $X$ and $Y$ switched.

The curvature tensor satisfies 
$$-\langle Z, R(X,Y)W\rangle = \langle R(X,Y)Z,W\rangle = \langle R(Z,W)X,Y
\rangle$$
 pointwise, so
\begin{eqnarray*} \langle A_XY,Z\rangle_1 &=&
  -\frac{1}{2}\langle[D_Z,\wgt]X,Y\rangle_0\\
&=& -\frac{1}{2}\langle (-\ndg(R(Z,\dg)X) - R(Z,\dg)\ndg X,Y\rangle_0\\
&=& -\frac{1}{2} \langle R(Z,\dg)X,\ndg Y\rangle_0 + \frac{1}{2}\langle
  R(Z,\dg)\ndg X,Y\rangle_0\\
&=& -\frac{1}{2} \langle R(X,\ndg Y)Z,\dg\rangle_0 + \frac{1}{2}\langle R(\ndg X,
  Y)Z,\dg\rangle_0\\
&=& \frac{1}{2}\langle Z, R(X,\ndg Y)\dg\rangle_0 - \frac {1}{2} \langle Z,
  R(\ndg X,Y)\dg\rangle_0\\
&=&\frac{1}{2}\langle Z, \wgti(R(X,\ndg Y)\dg - R(\ndg X, Y)\dg)\rangle_1.
\end{eqnarray*}
Thus $A_XY$ must equal $\frac{1}{2} \wgti(R(X,\ndg Y)\dg - R(\ndg X, Y)\dg).$  
This makes sense: for $X, Y\in H^{s'-1}$, 
 $A_XY\in H^{s'}\subset H^1,$ since $R$ is zeroth order.  
\hfill$\Box$
\medskip

\begin{rem} Locally on $LM$, we should have 
$D_XY = \delta_X^{LM}Y +  \omega_X^{LM}(Y)$.  
Now
$\delta_X^{LM}Y$  can only mean $\frac{d}{d\tau}|_{\tau =
  0}\frac{d}{d\epsilon}|_{\epsilon = 0}\gamma(\epsilon,\tau,\theta)$, where
$\gamma(0,0,\theta) = \gamma(\theta)$, ${d\epsilon}|_{\epsilon =
  0}\gamma(\epsilon,0,\theta) = X(\theta)$, 
${d\tau}|_{\tau = 0}\gamma(\epsilon,\tau,\theta) = Y_{\gamma(\epsilon, 0,\cdot)}(\theta).$
 In other words, $\delta_X^{LM}Y$ equals $ \delta_XY$.
Since $D_XY^a = \delta_XY^a + \chw{b}{c}{a}(\gamma(\theta))$, the connection one-form for the $L^2$ Levi-Civita connection on $LM$ is given by
$$\omega^{LM}_X(Y)^a(\gamma)(\theta) = \chw{b}{c}{a}(\gamma(\theta))X^bY^c
= \omega^M_X(Y)^a(\gamma(\theta)).$$
\end{rem}

By this remark, we get
\begin{cor}\label{cor2} The connection one-form $\omega = \omega^1$ for $\nabla^1$ in 
exponential coordinates is
\begin{eqnarray}\label{two}\omega_X(Y)(\gamma)(\theta) &=& \omega^M_X(Y)(\gamma(\theta))
  +  
\frac{1}{2}\bigl\{\wgti\left[
-\ndg(R(X,\dg)Y) - R(X,\dg)\ndg Y\right.\nonumber\\
&&\qquad \left. -\ndg(R(Y,\dg)X) - R(Y,\dg)\ndg X\right.\\
&&\qquad  \left.+R(X,\ndg Y)\dg - R(\ndg X, Y)\dg\right]\bigr\}(\theta).\nonumber
\end{eqnarray}
\end{cor}

\subsection{The Levi-Civita Connection for $s\in\BbZ^+$}
${}$
\medskip

For $s>\frac{1}{2}$, the proof of Prop.~\ref{old1.3} extends directly to give

\begin{lem} \label{lem: LCs}
The Levi-Civita connection for the $H^s$ metric is given by
$$\nabla_X^sY = D_XY + \frac{1}{2}\wgtsi[D_X, \wgts]Y +
  \frac{1}{2}\wgtsi[D_Y, \wgts]X 
+ A_XY,$$
where we assume that for $X, Y\in H^{s'-1}$, $A_XY\in H^s$ is characterized by
\begin{equation}\label{axy}
-\frac{1}{2}\langle [D_Z,\wgts]X,Y\rangle_0 = \langle A_XY,Z\rangle_s.
\end{equation}
\end{lem}
\bigskip

We now compute the bracket terms.

\begin{lem}\label{bracketterms}
For $s\in \Z^+$, at the loop $\gamma$,
\begin{equation}\label{bracket}
[D_X,\wgts]Y = \sum_{k=1}^s(-1)^k\left(\begin{array}{c}s\\k\end{array}\right)
\sum_{j=0}^{2k-1} \nabla_{\dg}^j(R(X,\dg)\nabla_{\dg}^{2k-1-j}Y).
\end{equation}
In particular, $[D_X,\wgts]Y$ is a $\pdo$ of order $2s-1$ in either $X$ or $Y$.
\end{lem}

\begin{proof}  The sum over $k$ comes from the binomial expansion of $\wgts$, so
we just need an inductive formula for 
$[D_X,\Delta^s].$   
The case $s=1$ is Proposition \ref{old1.3}. For the induction step, we have
\begin{eqnarray*} [D_X,\Delta^s] &=& D_X\Delta^{s-1}\Delta - \Delta^sD_X\\
&=& \Delta^{s-1}D_X\Delta + [D_X,\Delta^{s-1}]\Delta - \Delta^sD_X\\
&=& \Delta^sD_X +\Delta^{s-1}[D_X,\Delta] + [D_X,\Delta^{s-1}]\Delta
-\Delta^sD_X\\
&=& \Delta^{s-1}(-\nabla_{\dg}(R(X,\dg)Y) -R(X,\dg)\nabla_{\dg}Y)\\
&&\qquad  -  
\sum_{j=0}^{2s-3}(-1)^{s-1}
\nabla^j_{\dg}(R(X,\dg)\nabla_{\dg}^{2k-j-1}(-\nabla^2_{\dg}Y)\\
&=& (-1)^{s-1}(-\nabla_{\dg}^{2s-1}(R(X,\dg)Y) - (-1)^{s-1}\nabla_{\dg}^{2s-2}(R(X,\dg)\nabla_{\dg}Y)\\
&&\qquad + \sum_{j=0}^{2s-3}(-1)^{s}
\nabla^j_{\dg}(R(X,\dg)\nabla_{\dg}^{2k-j-1}(-\nabla^2_{\dg}Y)\\
&=& \sum_{j=0}^{2s-1}(-1)^s \nabla_{\dg}^j(R(X,\dg)\nabla_{\dg}^{2k-1-j}Y).
\end{eqnarray*} 
\end{proof}

We check that $A_XY$ is a $\pdo$ in $X$ and $Y$ for $s\in \BbZ^+.$

\begin{lem} \label{insert3} For $s\in\BbZ^+$ and fixed $X, Y\in H^{s'-1}$, $A_XY$ in (\ref{axy})
 is an explicit $\pdo$ in $X$ and $Y$ of order at most $-1.$
\end{lem}

\begin{proof}  By (\ref{bracket}), for  $j, 2k-1-j \in \{0,1,...,2s-1\}$, a typical term on
the left hand side of (\ref{axy}) is
\begin{eqnarray*}   \ipo{\nabla^j_{\dg}(R(Z,\dg)\nabla_{\dg}^{2k-1-j}X)}{Y} &=& 
(-1)^j
 \ipo{R(Z,\dg)\nabla_{\dg}^{2k-1-j}X}{\nabla^j_{\dg} Y}\\
&=& (-1)^j\ints g_{i\ell} (R(Z,\dg)\nabla_{\dg}^{2k-1-j}X)^i(\nabla^j_{\dg} Y)^\ell d\theta\\
&=& (-1)^j\ints g_{i\ell} Z^k R_{krn}^{\ \ \ i}\dg^r (\ndg^{2k-1-j}X)^n (\ndg^jY)^\ell d\theta\\
&=& (-1)^j
 \ints g_{tm}g^{kt}  g_{i\ell} Z^m R_{krn}^{\ \ \ i}\dg^r (\ndg^{2k-1-j}X)^n (\ndg^jY)^\ell d\theta\\
&=& (-1)^j \ipo{Z}
{g^{kt} g_{i\ell} R_{krn}^{\ \ \ i}\dg^r (\ndg^{2k-1-j}X)^n (\ndg^jY)^\ell\pa_t}\\ 
&=& (-1)^j\ipo{Z}{R^t_{\ rn\ell} \dg^r (\ndg^{2k-1-j}X)^n (\ndg^jY)^\ell\pa_t}\\ 
&=&(-1)^{j+1} \ipo{Z}{R^{\ \ \ t}_{n\ell r} \dg^r (\ndg^{2k-1-j}X)^n (\ndg^jY)^\ell\pa_t}\\
&=& (-1)^{j+1}\ipo{Z}{R(\ndg^{2k-1-j}X,\ndg^jY)\dg}\\
&=& (-1)^{j+1} \ips{Z}{\wgtsi R(\ndg^{2k-1-j}X,\ndg^jY)\dg}.
  \end{eqnarray*}
  (In the integrals and inner products, the local expressions are in fact globally defined one-forms on $S^1$, resp.~vector fields along $\gamma$, so we do not need a partition of unity.)
$\wgtsi R(\ndg^{2k-1-j}X,\ndg^jY)\dg$ is of order at most $-1$ in either $X$ or $Y$, so this term is in 
$H^{s'}\subset H^s.$  Thus the last inner product is well defined.  
\end{proof}

By  (\ref{axy}), (\ref{bracket}) and the proof of Lemma \ref{insert3}, we get
$$A_XY = \sum_{k=1}^s(-1)^k\left(\begin{array}{c}s\\k\end{array}\right)
\sum_{j=0}^{2k-1} (-1)^{j+1}   \wgtsi R(\ndg^{2k-1-j}X,\ndg^jY)\dg.$$
This gives:

\begin{thm} \label{thm:sinz}
For $s\in\BbZ^+$, the Levi-Civita connection for the $H^s$ metric at the
loop $\gamma$ is given by
\begin{eqnarray*} \nabla_X^sY(\gamma) &=& D_XY(\gamma) + \frac{1}{2}\wgtsi
 \sum_{k=1}^s(-1)^k\left(\begin{array}{c}s\\k\end{array}\right)
\sum_{j=0}^{2k-1} \nabla_{\dg}^j(R(X,\dg)\nabla_{\dg}^{2k-1-j}Y)\\
&&\qquad + \xly\\
  &&\qquad 
 +  \sum_{k=1}^s(-1)^k\left(\begin{array}{c}s\\k\end{array}\right)
\sum_{j=0}^{2k-1} (-1)^{j+1}   \wgtsi R(\ndg^{2k-1-j}X,\ndg^jY)\dg.
\end{eqnarray*}
\end{thm}

\subsection{The Levi-Civita Connection for General $s>\frac{1}{2}$}

${}$
\medskip

In this subsection, we show that the $H^s$ Levi-Civita connection for general $s>\frac{1}{2}$ exists in the strong sense of Remark \ref{lcrem}.
The formula is less explicit than in the 
$s\in \Z^+$ case, but is good enough for symbol calculations.

By Lemma \ref{lem: LCs}, we have to examine the term $A_XY$, which, if it exists, is
 characterized by (\ref{axy}):
$$-\frac{1}{2}\ipo{[D_Z,\weight]X}{Y} = \ips{A_XY}{Z}$$
for $Z\in H^s$.  As explained in Remark \ref{lcrem}, we may take
$X, Y\in H^{s'-1}.$
Throughout this section we assume that $s'\gg s$.

The following lemma extends Lemma \ref{bracketterms}.
\begin{lem}\label{pdo}
 (i)   For fixed $Z\in H^{s'-1}$,  $[D_Z,\weight] X$ is a $\Psi$DO of
  order $2s-1$ in $X$. For ${\rm Re}(s)\neq 0$, the principal symbol of $[D_Z,\weight]$ is 
  linear in $s$.
  
  (ii) For fixed $X\in H^{s'-1}$, $[D_Z,\weight]X$ is a $\pdo$ 
  of order $2s-1$   in $Z$.
\end{lem}

As usual, ``of order $2s-1$" means ``of order at most $2s-1.$"

\begin{proof}
(i) For $f:LM\to \BbC$, we get $[D_Z,\weight]fX = f[D_Z,\weight]X$, since $[f,\weight]=0.$  
Therefore, $[D_Z,\weight]X$ depends only on $X|_\gamma.$

By Lemma \ref{lem:l2lc}, $D_Z = \delta_Z + \Gamma \cdot Z$ in shorthand exponential 
coordinates.  The Christoffel symbol term is zeroth order and $\weight$ has scalar leading order symbol, so $[\Gamma\cdot Z,\weight]$ has order $2s-1.$  

From the integral expression for 
$\weight$, it is immediate that 
\begin{eqnarray}\label{immediate}
[\delta_Z,\weight]X &=& (\delta_Z\weight) X + \weight\delta_Z X - \weight\delta_ZX\\
&=& (\delta_Z\weight) X.\nonumber
\end{eqnarray}
$\delta_Z\weight$ is a limit of differences of $\pdo$s on bundles isomorphic to $\gamma^*TM$.
Since the algebra of $\pdo$s is closed in the Fr\'echet topology
of all $C^k$ seminorms
of symbols and smoothing terms
on compact sets, $\delta_Z\weight$ is a $\pdo.$

Since $\weight$ has order $2s$ and has scalar leading order symbol, 
$[D_Z,\weight]$ have order $2s-1$.  For later purposes (\S3.2), we compute some explicit symbols.  

Assume Re$(s)<0.$  As in the construction of $\weight$,
 we will compute what the symbol asymptotics
of $\delta_Z\weight$ should
be, and then construct an operator with these asymptotics.
From the functional calculus for unbounded operators, we have
\begin{eqnarray}\label{funcalc}
\delta_Z\weight &=& \delta_Z\left(\frac{i}{2\pi}\int_\Gamma
\lambda^s(1+\Delta-\lambda)^{-1}d\lambda\right)\nonumber\\
&=& \frac{i}{2\pi}\int_\Gamma
\lambda^s\delta_Z (1+\Delta-\lambda)^{-1}d\lambda\\
&=& -\frac{i}{2\pi}\int_\Gamma
\lambda^s (1+\Delta-\lambda)^{-1} (\delta_Z\Delta)
(1+\Delta-\lambda)^{-1}d\lambda,\nonumber
\end{eqnarray}
where $\Gamma$ is a contour around the spectrum of $1+\Delta$, and the
hypothesis on $s$ justifies the exchange of $\delta_Z$ and the integral.  The
operator $A =
(1+\Delta-\lambda)^{-1} \delta_Z\Delta  (1+\Delta-\lambda)^{-1}$ is a $\pdo$
of order $-3$
 with top order symbol
\begin{eqnarray*} \sigma_{-3}(A)(\theta,\xi)^\ell_j &=&
(\xi^2-\lambda)^{-1}\delta^\ell_k (-2Z^i\partial_i\ch{\nu}{\mu}{k}\dg^\nu
   -2\ch{\nu}{\mu}{k}\dot Z^\nu) \xi (\xi^2-\lambda)^{-1}\delta^\mu_j\\
&=&
(-2Z^i\partial_i\ch{\nu}{j}{\ell}\dg^\nu
   -2\ch{\nu}{j}{\ell}\dot Z^\nu)
\xi (\xi^2-\lambda)^{-2}.
\end{eqnarray*}
Thus the top order symbol of $\delta_Z\weight$ should be
\begin{eqnarray}\label{tsmo}
 \sigma_{2s-1}(\delta_Z\weight)(\theta,\xi)^\ell_j
&=& -\frac{i}{2\pi}\int_\Gamma
\lambda^s (-2Z^i\partial_i\ch{\nu}{j}{\ell}\dg^\nu
   -2\ch{\nu}{j}{\ell}\dot Z^\nu)
\xi (\xi^2-\lambda)^{-2} d\lambda  \nonumber\\
&=& \frac{i}{2\pi}\int_\Gamma s\lambda^{s-1}
(-2Z^i\partial_i\ch{\nu}{j}{\ell}\dg^\nu
   -2\ch{\nu}{j}{\ell}\dot Z^\nu)
\xi (\xi^2-\lambda)^{-1} d\lambda \nonumber\\
&=& s(-2Z^i\partial_i\ch{\nu}{j}{\ell}\dg^\nu
   -2\ch{\nu}{j}{\ell}\dot Z^\nu)\xi (\xi^2-\lambda)^{s-1}.
\end{eqnarray}
Similarly, all the terms in the symbol asymptotics for $A$ are of the form
$B^\ell_j \xi^n(\xi^2-\lambda)^m$ for some matrices $B^\ell_j =
B^\ell_j(n,m).$   This produces a symbol sequence $
\sum_{k\in \Z^+}\sigma_{2s-k}$,  and there exists a  $\pdo$ $P$ with $\sigma(P) =
\sum \sigma_{2s-k}$.  (As in \S\ref{pdoreview}, we 
first produce operators $P_i$ on a coordinate cover $U_i$ of $S^1$, 
and then set
$P = \sum_i\phi_iP_i\psi_i$.) The construction
depends on
the choice of local coordinates
covering $\gamma$, the partition of unity and cutoff
functions as above, and a cutoff function in $\xi$; as
usual, different choices change the operator by a smoothing operator.
Standard estimates 
show that $P-\delta_Z\weight$ is a smoothing
operator, this verifies explicitly that  $\delta_Z\weight$ is a $\pdo$ of order $2s-1.$

For Re$(s) >0$,  motivated by differentiating $\weightinv\circ\weight = {\rm
  Id}$, we set
\begin{equation}\label{abc}
\delta_Z\weight = -\weight\circ\delta_Z\weightinv\circ\weight.
\end{equation}
This is again a $\pdo$ of order $2s-1$ with principal symbol linear in $s$. 

 (ii) As a $\pdo$ of order $2s$, $\weight$ has the expression
 $$\weight X(\gamma)(\theta) = \int_{T^*S^1}
 e^{i(\theta-\theta')\cdot \xi} 
 p(\theta,\xi) X(\gamma)(\theta')d\theta' d\xi,$$
 where we omit the cover of $S^1$ and its partition of unity on the right hand side.
 Here $
 p(\theta,\xi)$ is the symbol of $\weight$, which has the asymptotic expansion
 $$
 p(\theta,\xi) \sim \sum_{k=0}^\infty p_{2s-k }(\theta,\xi).$$
 The covariant derivative along $\gamma$ on
$Y\in\Gamma(\gamma^*TM)$ is given by
\begin{eqnarray*}\frac{DY}{d\gamma} &=&
(\gamma^*\nabla^{M})_{\partial_\theta}(Y) =
\partial_\theta Y + (\gamma^*\omega^{M})(\partial_\theta)(Y)\\
&=& \partial_\theta(Y^i)\partial_i + \dot\gamma^t Y^r
\Gamma^j_{tr}\partial_j,
\end{eqnarray*}
where $\nabla^{M}$ is the Levi-Civita connection on $M$ and $\omega^{M}$ is the
connection one-form in exponential coordinates on $M$. 
For $\Delta =
(\frac{D}{d\gamma})^* \frac{D}{d\gamma}$, an integration by parts using the
formula
$\partial_tg_{ar} = \Gamma_{\ell t}^ng_{rn} + \Gamma_{rt}^ng_{\ell n}$  gives
$$(\Delta Y)^k = -\partial^2_\theta Y^k
-2\Gamma_{\nu\mu}^k\dot\gamma^\nu\partial_\theta Y^\mu -\left
( \partial_\theta\Gamma_{\nu\delta}^k\dot\gamma^\nu
+\Gamma_{\nu\delta}^k\ddot\gamma^\nu +
\Gamma_{\nu\mu}^k\Gamma_{\e\delta}^\mu\dot\gamma^\e\dot\gamma^\nu\right)
Y^\delta.$$
Thus $p_{2s}(\theta, \xi) =  |\xi|^2$ is independent of $\gamma$, but the lower order symbols  depend on 
 derivatives of both $\gamma$ and the metric on $M$. 
 
 We have
 \begin{eqnarray} [D_Z,\weight]X(\gamma)(\theta) &=&
 D_Z \int_{T^*S^1}
 e^{i(\theta-\theta')\cdot \xi} 
 p(\theta,\xi) X(\gamma)(\theta')d\theta' d\xi\label{216}\\
 &&\quad - \int_{T^*S^1}
 e^{i(\theta-\theta')\cdot \xi} 
 p(\theta,\xi) D_ZX(\gamma)(\theta')d\theta' d\xi. \label{217}
 \end{eqnarray}
 In local coordinates, (\ref{216}) equals
 \begin{eqnarray} \label{218}
 \lefteqn{
\left[ D_Z \int_{T^*S^1}
 e^{i(\theta-\theta')\cdot \xi} 
 p(\theta,\xi) X(\gamma)(\theta')d\theta' d\xi\right]^a}\nonumber\\
 &=& \delta_Z\left[
 \int_{T^*S^1}
 e^{i(\theta-\theta')\cdot \xi} 
 p(\theta,\xi) X(\gamma)(\theta')d\theta' d\xi\right]^a(\theta)\\
 &&\quad + \Gamma^a_{bc} Z^b(\gamma)(\theta) 
 \left[ \int_{T^*S^1}
 e^{i(\theta-\theta')\cdot \xi} 
 p(\theta,\xi) X(\gamma)(\theta')d\theta' d\xi\right]^c(\theta).\nonumber
 \end{eqnarray}
 Here we have suppressed matrix indices in $p$ and $X$.
We can bring $\delta_Z$ past the integral on the right hand side of (\ref{218}).  If
$\gamma_\epsilon$ is a family of curves with $\gamma_0 = \gamma, \dot\gamma_\epsilon = Z$, then
$$\delta_Zp(\theta, \xi) = \frac{d}{d\epsilon}\biggl|_{_{_{\epsilon=0}}} p(\gamma_\epsilon,
\theta,\xi) = \frac{d\gamma_\epsilon^k}{d\epsilon}\biggl|_{_{_{\epsilon=0}}}
\partial_k p(\gamma,\theta, \xi) = Z^k(\gamma(\theta))\cdot \partial_k p(\theta,\xi).$$ Substituting this into (\ref{218}) gives
\begin{eqnarray}\label{219}
\lefteqn{
\left[ D_Z \int_{T^*S^1}
 e^{i(\theta-\theta')\cdot \xi} 
 p(\theta,\xi) X(\gamma)(\theta')d\theta' d\xi\right]^a}\nonumber\\
&=& 
\lefteqn{
\left[ \int_{T^*S^1}
 e^{i(\theta-\theta')\cdot \xi} 
 Z^k(\gamma)(\theta)\cdot \partial_k p(\theta,\xi) X(\gamma)(\theta')d\theta' d\xi\right]^a}\\
 &&\quad 
+ \Gamma^a_{bc} Z^b(\gamma)(\theta) 
 \left[ \int_{T^*S^1}
 e^{i(\theta-\theta')\cdot \xi} 
 p(\theta,\xi) X(\gamma)(\theta')d\theta' d\xi\right]^c(\theta).\nonumber\\
&&\qquad + 
 \left[ \int_{T^*S^1}
 e^{i(\theta-\theta')\cdot \xi} 
 p(\theta,\xi) \delta_Z X(\gamma)(\theta')d\theta' d\xi\right]^c(\theta).\nonumber
\end{eqnarray}
Similarly, (\ref{217}) equals
\begin{eqnarray}\label{220}
\lefteqn{\left[\int_{T^*S^1}
 e^{i(\theta-\theta')\cdot \xi} 
 p(\theta,\xi) D_ZX(\gamma)(\theta')d\theta' d\xi\right]^a}\nonumber\\
 &=& \left[\int_{T^*S^1}
 e^{i(\theta-\theta')\cdot \xi} 
 p(\theta,\xi) \delta_ZX(\gamma)(\theta')d\theta' d\xi\right]^a\\
 &&\quad + 
 \int_{T^*S^1}
 e^{i(\theta-\theta')\cdot \xi} 
 p(\theta,\xi)^a_e\Gamma^e_{bc}Z^b(\gamma)(\theta') X^c(\gamma)(\theta')d\theta' d\xi.
 \nonumber
 \end{eqnarray}
 Substituting (\ref{219}), (\ref{220}), into (\ref{216}), (\ref{217}), respectively, gives
 \begin{eqnarray}\label{221}
\lefteqn{( [D_Z,\weight]X(\theta))^a}\\
&=&
Z^b(\theta)\cdot\left[\int_{T^*S^1} e^{i(\theta-\theta')\cdot\xi}\left(\partial_bp^a_e(\theta,\xi) +
\Gamma_{bc}^a(\gamma(\theta)p_e^c(\theta, \xi)\right)X^e(\theta')d\theta'd\xi\right]\nonumber\\
&&\quad - \int_{T^*S^1} e^{i(\theta-\theta')\cdot\xi}p(\theta,\xi)^a_e\Gamma_{bc}^e
(\gamma(\theta'))Z^b(\theta') X^c(\theta')d\theta' d\xi,\nonumber
  \end{eqnarray}
 where $X(\theta') = X(\gamma)(\theta)$ and similarly for Z.
 
 The first term on the right hand side of (\ref{221}) is order zero in $Z$; note that
 $0<2s-1$, since $s>\frac{1}{2}$.  For the last term in (\ref{221}), we do a change of variables typically used in the proof that the composition of $\pdo$s is a $\pdo.$  Set
 \begin{equation}\label{221a}q(\theta, \theta', \xi)^a_b = p(\theta,\xi)^a_e \Gamma_{bc}^e(\gamma(\theta'))X^c
 (\theta'),
 \end{equation}
 so the last term equals
 \begin{eqnarray*}
(PZ)^a(\theta) &\stackrel{\rm def}{=}&  \int_{T^*S^1} e^{i(\theta-\theta')\cdot\xi}q(\theta, \theta', \xi)^a_b Z^b(\theta') d\theta' d\xi\\
 &=& \int_{T^*S^1} e^{i(\theta-\theta')\cdot\xi}q(\theta, \theta', \xi)^a_b e^{i(\theta'-\theta'')
 \cdot\eta} Z^b(\theta'') d\theta'' d\eta \ d\theta' d\xi,
 \end{eqnarray*}
 by applying  Fourier transform and its inverse to $Z$. A little algebra gives
 \begin{equation}\label{222}
 (PZ)^a(\theta) = \int_{T^*S^1} e^{i(\theta-\theta')\cdot\eta}r(\theta,\eta)^a_b Z^b(\theta')
 d\theta' d\eta,
 \end{equation}
 with 
 \begin{eqnarray*}r(\theta, \eta) &=& \int_{T^*S^1} e^{i(\theta-\theta')\cdot(\xi-\eta)}
 q(\theta,\theta', \xi) d\theta' d\xi\\
 &=&  \int_{T^*S^1} e^{it\cdot\xi} q(\theta,\theta - t, \eta + \xi) dt\  d\xi.
 \end{eqnarray*} 
 In the last line we continue to abuse notation by treating the integral in local coordinates in 
 $t = \theta-\theta'$ lying in an interval $I\subset \R$ and implicitly
 summing over a cover and partition of unity of $S^1;$ thus we can consider $q$ as a compactly supported function in $t\in\R.$
 Substituting in the Taylor expansion of $q(\theta,\theta - t, \eta + \xi)$ in $\xi$ gives in local coordinates
  \begin{eqnarray}\label{223a}
  r(\theta, \eta) &=& \int_{T^*\R} e^{it\cdot \xi} \left[ \sum_{\alpha, |\alpha|=0}^N
 \frac{1}{\alpha!} \partial_\xi^\alpha|_{\xi=0} q(\theta,\theta-t, \eta+\xi)\xi^\alpha + {\rm O}
 (|\xi|^{N+1})\right] dt \  d\xi\nonumber\\
 &=& \sum_{\alpha, |\alpha|=0}^N \frac{i^{|\alpha|}}{\alpha!} \partial^\alpha_t\partial^\alpha
 _\xi q(\theta, \theta, \eta) + {\rm O} (|\xi|^{N+1}).
 \end{eqnarray}
 Thus $P$ in (\ref{222}) is a $\pdo$ with apparent top order symbol 
 $q(\theta, \theta, \eta)$, which by (\ref{221a}) has order $2s.$  The top order symbol can be computed in any local coordinates on $S^1$ and $\gamma^*TM$.  If we choose
 manifold coordinates (see \S2.3) which are 
Riemannian normal coordinates centered at $\gamma(\theta)$, the Christoffel symbols vanish at this point, and
 so
 $$q(\theta, \theta, \eta)^a_b = p(\theta,\xi)^a_e\Gamma_{bc}^e(\gamma(\theta)) X^c(\theta)
 =0.$$
 Thus $P$ is in fact of order $2s-1$, and so both terms on the right hand side of (\ref{221}) have order at most $2s-1$.

\end{proof}

\begin{rem} (i) For $s\in\Z^+$,  $\delta_Z\weight$ differs
from the usual definition by a smoothing operator.  

(ii) For all $s$, the proof of Lemma \ref{pdo}(i) shows that
$\sigma(\delta_Z\weight) = \delta_Z(\sigma(\weight)).$

\end{rem}

We can now complete the computation of the Levi-Civita connection for general $s.$

Let $[D_\cdot,\weight]X^*$ be the formal $L^2$ adjoint  of $[D_\cdot,\weight]X$.
We abbreviate $[D_\cdot,\weight]X^*(Y)$ by $[D_Y,\weight]X^*.$

\begin{thm}  \label{thm25} (i) For $s>\frac{1}{2}$, 
The Levi-Civita connection for the $H^s$ metric is given by
\begin{eqnarray}\label{quick}\nabla_X^sY &=& D_XY + \frac{1}{2}\wgtsi[D_X, \wgts]Y +
  \frac{1}{2}\wgtsi[D_Y, \wgts]X \nonumber\\
&&\quad -\frac{1}{2} \wgtsi[D_Y,\weight]X^*.
\end{eqnarray}

(ii) The connection one-form $\omega^s$ in exponential coordinates is given by
\begin{eqnarray}\label{223}\lefteqn{\omega^s_X(Y)(\gamma) (\theta)}\\
&=& \omega^M(Y)(\gamma(\theta)) + 
\left(\frac{1}{2}\wgtsi[D_X, \wgts]Y +
  \frac{1}{2}\wgtsi[D_Y, \wgts]X \right.\nonumber\\
  &&\quad \left.
-\frac{1}{2} \wgtsi[D_Y,\weight]X^*\right)(\gamma)(\theta).\nonumber
\end{eqnarray}

(iii) The connection one-form takes values in zeroth order $\pdo$s.
\end{thm}

\begin{proof}  Since $[D_Z,\weight]X$ is a $\pdo$ in $Z$ of order $2s-1$, its formal adjoint is
a $\pdo$ of the same order.  Thus
$$\langle [D_Z,\weight]X,Y\rangle_0 = \langle Z, [D_\cdot, \weight]X^*(Y)\rangle
= \langle Z, \wgtsi[D_Y,\weight]X^*\rangle_s.$$
Thus $A_XY$ in (\ref{axy}) satisfies
$A_XY = \wgtsi[D_Y,\weight]X^*.$  Lemma \ref{lem: LCs} applies to all $s>\frac{1}{2}$, 
so (i) follows.  (ii) follows as 
in Corollary \ref{cor2}.  Since $\omega^M$ is zeroth order and all 
other terms have order $-1$, (iii) holds as well.
\end{proof}

\begin{rem}  This theorem implies that the Levi-Civita connection exists for the 
$H^s$ metric in the strong sense:  for $X\in  T_\gamma LM =H^{s'-1}(\gamma^*TM)$
and $Y\in H^{s'-1}(\cdot^*TM)$ a smooth vector field on $LM = H^{s'}(S^1,M)$,
 $\nabla^s_XY(\gamma)\in H^{s'-1}(\gamma^*TM).$  (See Remark
2.1.)  For each term except $D_XY$ on the right hand side of (\ref{quick}) is order
$-1$ in $Y$, and so takes $H^{s'-1}$ to $H^{s'}\subset H^{s'-1}.$  For $D_XY = \delta_XY + \Gamma\cdot Y$, $\Gamma$ is zeroth order and so bounded on $H^{s'-1}.$  Finally, 
the definition of a smooth vector field on $LM$ implies that $\delta_XY$ stays in $H^{s'-1}$
for all $X$.
\end{rem}

\subsection{{\bf Extensions of the Frame Bundle of $LM$}}\label{extframe}

In this subsection we discuss the choice of structure group for the
$H^s$ and Levi-Civita connections on $LM.$

Let $\calH$ be the Hilbert space
 $H^{s_0}(\gamma^*TM)$ for 
a fixed $s_0$ and $\gamma.$ 
 Let $GL(\calH)$ be the group of bounded invertible linear
operators on $\calH$; inverses of elements are bounded by the closed graph
theorem.   $GL(\calH)$ has the subset
topology of the norm topology on ${\mathcal B}(\calH)$, the bounded linear
operators on $\calH$.
$GL(\calH)$ is an infinite dimensional Banach Lie group, as a group which
is an open subset of the infinite dimensional Hilbert manifold 
${\mathcal B}(\calH)$
\cite[p.~59]{Omori}, and has Lie algebra 
${\mathcal B}(\calH)$. Let $\pdo_{\leq 0}, 
\pdo_0^*$ denote the algebra of classical
$\pdo$s of nonpositive order 
and the group of invertible zeroth order $\pdo$s, respectively,
where all $\pdo$s act on $\calH.$   
Note that $\pdo_0^*\subset GL(\calH).$
   
\begin{rem} 
The inclusions of $\pdo_0^*, \pdo_{\leq 0}$ into $GL(\calH), {\mathcal
    B}(\calH)$ are trivially continuous in the subset topology.
For the Fr\'echet topology on $\pdo_{\leq 0}$, 
the  inclusion is 
continuous as in \cite{lrst}.
\end{rem}

We recall
the relationship between
 the  connection one-form $\theta$ on the frame bundle $FN$ of a
 manifold $N$
and
local expressions for the connection on $TN.$ For $U\subset N$,
 let $\chi:U\to FN$ be a local section.
 A metric  connection $\nabla$ on $TN$ with local
connection one-form $\omega$ determines a connection $\theta_{FN}\in
 \Lambda^1(FN, {\mathfrak o}(n))$ on $FN$
by {\it (i)} $\theta_{FN}$ is the Maurer-Cartan one-form on each fiber,
and {\it (ii) }
$\theta_{FN}(Y_u)=\omega (X_p),$ for $ Y_u=\chi_*X_p$
\cite[Ch.~8, Vol.~II]{Spi}, or equivalently
$\chi^*\theta_{FN} = \omega.$

This applies to $N=LM.$
The frame bundle $FLM\to LM$ is constructed
as in the finite dimensional case. The
fiber over $\gamma$ is isomorphic to the gauge group $\calG$ of $\calR$
and fibers are glued by the transition functions for
$TLM$. Thus the frame bundle is
topologically a
$\calG$-bundle.

However, by Theorem \ref{thm25},
the Levi-Civita connection one-form $\omega^s_X$
takes
values in $\pdo_{\leq 0}$. 
The curvature two-form $\Omega^{s} = d_{LM}\omega^{s} + \omega^{s}\wedge
\omega^s$ also takes values in $\pdo_{\leq 0}.$  (Here $d_{LM}\omega^{s}(X,Y)$
is defined by the Cartan formula for the exterior derivative.)
These
forms should take values in the Lie algebra of the structure
group.  Thus we should extend the structure group to the Fr\'echet Lie group
 $\pdo_0^*$, since its Lie
algebra is $\pdo_{\leq 0}.$  
This leads to an extended
frame bundles, also denoted $FLM$. The  transition
 functions are unchanged, since 
$\calG \subset \pdo_0^*$.
 Thus $(FLM,\theta^s)$ as a geometric
bundle (i.e.~as a  bundle with connection $\theta^s$ associated to
$\nabla^{1,s}$) is a $\pdo_0^*$-bundle.

In summary, for the Levi-Civita connections we have
$$ \begin{array}{ccc}
\calG&\longrightarrow &FLM\\
& & \downarrow\\
& & LM
\end{array}
\ \ \ \ \ \ \ \ \ \ \ \ \ \ 
\begin{array}{ccc}
\pdo_0^*&\longrightarrow &(FLM,\theta^s)\\
& & \downarrow\\
& & LM
\end{array}
$$

\begin{rem}\label{rem:ext}  If
 we extend the structure group of the frame bundle with
  connection from $\pdo_0^*$ to
  $GL(\calH)$, the frame bundle becomes trivial by Kuiper's theorem.  
Thus
there is a potential loss of information if 
we pass to the larger frame
  bundle.

The situation is similar to the following examples.  Let $E\to S^1$ be
the $GL(1,\R)$ (real line)
 bundle with gluing functions (multiplication by) $1$ at $1\in
S^1$ and $2$ at $-1\in S^1.$  $E$ is trivial as a $GL(1,\R)$-bundle, 
with global section $f$ with $\lim_{\theta\to -\pi^+}f(e^{i\theta}) = 1, 
f(1) = 1,
\lim_{\theta\to\pi^-}f(e^{i\theta}) = 1/2.$  
However, as a $GL(1,\Q)^+$-bundle, $E$ is nontrivial, as a
global section is locally constant. As a second example,
 let $E\to M$ be a nontrivial
$GL(n,\C)$-bundle. Embed $\C^n$ into a Hilbert space $\calH$, and extend $E$
to an $GL(\calH)$-bundle $\calE$ 
with fiber $\calH$ and with the  transition functions for $E$ (extended by the identity in
directions perpendicular to the image of $E$).  Then $\calE$ is
trivial.

\end{rem}

\section{{\bf Local Symbol Calculations}}\label{localsymbols}

In this section, we  compute the $0$ and $-1$ order symbols of the
connection one-form and the curvature two-form of
 the $s=1$ Levi-Civita connection. 
 We also compute the $0$ and $-1$ order symbols of the
connection one-form for the general $s>\frac{1}{2}$ connection, and the $0$ order symbol of the 
curvature of the general $s$ connection.
 These results are used in the calculations of Wodzicki-Chern-Simons
classes in \S6. The formulas  show that
the $s$-dependence of these symbols is
linear, which will be used to define regularized Wodzicki-Chern-Simons classes
(see Definition \ref{def:regularized}).

\subsection{{\bf Connection and Curvature Symbols for $s=1$}}
${}$
\medskip

In this subsection $\omega = \omega^1, \Omega = \Omega^1.$

Using Corollary \ref{cor2}, we can compute these symbols easily. 

\begin{lem} \label{old2.1}
(i) At $\gamma(\theta)$,
$\sigma_0(\omega_X)^a_b =  (\omega^M_X)^a_b = \chw{c}{b}{a}X^c.$

(ii)  \begin{eqnarray*}
\frac{1}{i|\xi|^{-2}\xi}\sigma_{-1}(\omega_X) &=& \frac{1}{2}(-2R(X,\dg)
-R(\cdot,\dg)X + R(X,\cdot)\dg).
\end{eqnarray*}
Equivalently,
\begin{eqnarray*}
\frac{1}{i|\xi|^{-2}\xi}\sigma_{-1}(\omega_X)^a_b &=& 
\frac{1}{2}
(-2R_{cdb}^{\ \ \ a} -R_{bdc}^{\ \ \ a} + R_{cbd}^{\ \ \  a})X^c\dg^d.
\end{eqnarray*}
\end{lem}

\begin{proof} (i) For $\sigma_0(\omega_X)$, the only term in (\ref{two}) of order
  zero is the Christoffel term.  

(ii) For $\sigma_{-1}(\omega_X)$, label the last six terms on the right hand side of
(\ref{two}) by (a), ..., (f).  By Leibniz rule for the tensors, the only
terms of order $-1$ come from:
in (a), $-\nabla_{\dg}(R(X,\dg)Y) = -R(X, \dg) \nabla_{\dg}Y +$ lower order in
$Y$;
in (b), the term  $-R(X, \dg) \nabla_{\dg}Y$;
in (c), the term $-R(\nabla_{\dg}Y, \dg)X$;
in (e), the term $R(X,\nabla_{\dg}Y)\dg.$

For any vectors $Z, W$, the curvature endomorphism $R(Z, W): TM\to TM$ has
$$R(Z,W)^a_b = R_{c d b}^{\ \ \ a}Z^cW^d.$$
  Also, since $(\nabla_{\dg}Y)^a =
\frac{d}{d\theta}Y^a $ plus zeroth order terms,
$\sigma_{1}(\nabla_{\dg}) =
i\xi\cdot Id.$
Thus in (a) and (b), 
$\sigma_1(-R(X, \dg) \nabla_{\dg})^a_b = -R_{cdb}^{\ \ \ a}X^c\dg^d\xi.$

For (c), we have $-R(\nabla_{\dg}Y, \dg)X = -R_{cdb}^{\ \ \
  a}(\nabla_{\dg}Y)^c\dg^d X^b\partial_a$, so the top order symbol is
$-R_{cdb}^{\ \ \ a}\xi\dg^dX^b = -R_{bdc}^{\ \ \ a}\xi\dg^d X^c.$

For (e), we have $R(X,\nabla_{\dg}Y)\dg = R_{cdb}^{\ \ \ a}X^c(\nabla_{\dg}Y)^d
\dg^b\partial_a$, so the top order symbol is 
$R_{cdb}^{\ \ \ a}X^c\xi \dg^b = R_{cbd}^{\ \ \ a}X^c\xi \dg^d.$

Since the top order symbol of $\wgti$ is $|\xi|^{-2}$, adding these four terms
finishes the proof. 
\end{proof}
 
 We now compute the top symbols of the curvature tensor.  $\sigma_{-1}(\Omega)$ involves
 the covariant derivative of the curvature tensor on $M$, but fortunately this symbol
 will not be needed in Part II.

\begin{lem}\label{old2.2}
(i) 
$\sigma_0(\Omega(X,Y))^a_b =  R^M(X,Y)^a_b = R_{cdb}^{\ \ \ a}X^cY^d.$

(ii) \begin{eqnarray*}
\frac{1}{i|\xi|^{-2}\xi}\sigma_{-1}(\Omega(X,Y)) &=&
\frac{1}{2}\left(\nabla_X[-2R(Y,\dg) - R(\cdot,\dg)Y +
  R(Y,\cdot)\dg]\right.\\
&&\qquad \left. - \xly \right.\\
&&\qquad\left. - [-2R([X,Y],\dg) -R(\cdot,\dg)[X,Y] + R([X,Y],\cdot)\dg] \right).
\end{eqnarray*}
Equivalently, in Riemannian normal coordinates on $M$ centered at $\gamma(\theta)$,
\begin{eqnarray}\label{moc}
\frac{1}{i|\xi|^{-2}\xi}\sigma_{-1}(\Omega(X,Y))^a_b &=& \frac{1}{2}
X[(-2R_{cdb}^{\ \ \ a} -R_{bdc}^{\ \ \ a} + R_{cbd}^{\ \ \
  a})\dg^d]Y^c - (X\leftrightarrow Y)\nonumber\\
&=&\frac{1}{2}
X[-2R_{cdb}^{\ \ \ a} -R_{bdc}^{\ \ \ a} + R_{cbd}^{\ \ \
  a}]\dg^dY^c  -(X\leftrightarrow Y)\\
&&\qquad +
\frac{1}{2}[-2R_{cdb}^{\ \ \ a} -R_{bdc}^{\ \ \ a} + R_{cbd}^{\ \ \
  a}]\dot X^dY^c - (X\leftrightarrow Y)\nonumber
\end{eqnarray}
\end{lem}

\begin{proof} 
(i)
\begin{eqnarray*} \sigma_0(\Omega(X,Y))^a_b &=& \sigma_0((d\omega +
  \omega\wedge\omega)(X,Y))^a_b\\
&=& [(d\sigma_0(\omega) + \sigma_0(\omega)\wedge\sigma_0(\omega))(X,Y)]^a_b\\
&=& [(d\omega^M + \omega^M\wedge\omega^M)(X,Y)]^a_b\\
&=& R^M(X,Y)^a_b = R_{cdb}^{\ \ \ a}X^cY^d.
\end{eqnarray*}

(ii) Since $\sigma_0(\omega_X)$ is independent of $\xi$, after dividing by
$i|\xi|^{-2}\xi$ we have
\begin{eqnarray*}\sigma_{-1}(\Omega(X,Y))^a_b &=& (d\sigma_{-1}(\omega)
  (X,Y))^a_b + \sigma_0(\omega_X)^a_c\sigma_{-1}(\omega_Y)^c_b
+ \sigma_{-1}(\omega_X)^a_c\sigma_{0}(\omega_Y)^c_b\\
&&\qquad
-\sigma_0(\omega_Y)^a_c\sigma_{-1}(\omega_X)^c_b
+ \sigma_{-1}(\omega_Y)^a_c\sigma_{0}(\omega_X)^c_b.
\end{eqnarray*}
As an operator on sections of $\gamma^*TM$, 
$\Omega^{LM} - \Omega^M$ has order $-1$  so $\sigma_{-1}(\Omega^{LM})
= \sigma_{-1}(\Omega^{LM} -\Omega^M)$ is independent of coordinates.
In Riemannian normal coordinates at $\gamma(\theta)$, $\sigma_0(\omega_X) = \sigma_0(\omega_Y) = 0$, so
\begin{eqnarray*}\sigma_{-1}(\Omega(X,Y))^a_b &=&
  (d\sigma_{-1}(\omega)(X,Y))^a_b\\
&=& X(\sigma_{-1}(\omega_Y))^a_b - Y(\sigma_{-1}(\omega_X))^a_b
  -\sigma_{-1}(\omega_{[X.Y]})^a_b\\
&=& \frac{1}{2} X[(-2R_{cdb}^{\ \ \ a} -R_{bdc}^{\ \ \ a} + R_{cbd}^{\ \ \
  a}]Y^c\dg^d] - (X\leftrightarrow Y)\\
&&\qquad -\frac{1}{2}( -2R_{cdb}^{\ \ \ a} -R_{bdc}^{\ \ \ a} + R_{cbd}^{\ \ \
  a}][X,Y]^c\dg^d.
\end{eqnarray*}
The terms involving $X(Y^c) - Y(X^c) - [X,Y]^c$ cancel (as they must, since the symbol two-form 
cannot involve derivatives of $X$ or $Y$).  Thus
$$\sigma_{-1}(\Omega(X,Y))^a_b = \frac{1}{2} X[(-2R_{cdb}^{\ \ \ a} 
-R_{bdc}^{\ \ \ a} + R_{cbd}^{\ \ \  a})Y^c\dg^d] - (X\leftrightarrow Y).$$

This gives the first coordinate expression in (\ref{moc}). The second expression follows from
 $X(\dg^d) = \dot X^d $ (see (\ref{badterms})).

To convert from the coordinate expression to the covariant expression, we follow the 
usual procedure of changing ordinary derivatives to covariant derivatives and adding bracket terms.  For example,
\begin{eqnarray*}\nabla_X(R(Y,\dg)) &=& (\nabla_XR)(Y,\dg) + R(\nabla_XY,\dg)
  + R(Y,\nabla_X\dg) \\
&=& X^iR_{cdb\ ;i}^{\ \ \ a}Y^c\dg^d + R(\nabla_XY,\dg) + R_{cdb}^{\ \ \
    a}Y^c(\nabla_X\dg)^d.
\end{eqnarray*}
In Riemannian normal coordinates at
$\gamma(\theta)$, we have $X^iR_{cdb\ ;i}^{\ \ \ a} = X^i\partial_i R_{cdb}^{\
  \ \ a} = X(R_{cdb}^{\ \ \ a})$ and $(\nabla_X\dg)^d =  X(\dg^d).$ 
  Thus 
  $$\nabla_X(R(Y,\dg)) -\xly - R([X,Y],\dg) = X(R_{cdb}^{\ \ \ a}\dg^d)Y^c - \xly.$$
  The other terms are handled similarly.
  \end{proof} 

\subsection{{\bf Connection and Curvature Symbols for General $s$}}
${}$
\medskip

The noteworthy feature of these computations is the linear dependence of $\sigma_{-1}(\omega^{s})$ on $s$.  

Let $g$ be the Riemannian metric on $M$.

\begin{lem}\label{lem:33}
(i) At $\gamma(\theta)$,
$\sigma_0(\omega^s_X)^a_b =  (\omega^M_X)^a_b = \chw{c}{b}{a}X^c.$

(ii) $\sigma_0(\Omega^s(X,Y))^a_b =  R^M(X,Y)^a_b = R_{cdb}^{\ \ \ a}X^cY^d.$

(iii)  $
\frac{1}{i|\xi|^{-2}\xi}\sigma_{-1}(\omega^s_X)^a_b = s T(X,\dg, g)$,
where $T(X, \dg, g)$  is tensorial and independent of $s$. 
\end{lem}

\begin{proof} (i)  By Lemma \ref{pdo}, the only term of order zero in (\ref{223}) 
 is $\omega^M_X.$

(ii)  The proof of Lemma \ref{old2.2}(ii) carries over.  

(iii) By Theorem \ref{thm25}, we have to compute $\sigma_{2s-1}$ for $[D_X,\wgts]$, 
$[D_\cdot,\wgts]X$, and $[D_\cdot,\wgts]X^*$, as 
$\sigma_{-1} (\wgtsi[D_X,\wgts]) =  |\xi|^{-2s}\sigma_{-1}([D_X,\wgts])$, etc.

Write $D_X = \delta_X + \Gamma\cdot X$ in shorthand.  Since $\wgts$ has scalar leading order symbol, $[\Gamma\cdot X,\wgts]$ has order $2s-1.$  Thus
we can compute $\sigma_{2s-1}([\Gamma\cdot X,\wgts])$ in any coordinate system.  
In  Riemannian normal coordinates centered at $\gamma(\theta)$, as in the proof of Lemma \ref{pdo}(ii),  the Christoffel symbols vanish.
Thus $\sigma_{2s-1}([\Gamma\cdot X,\wgts]) =0.$

By (\ref{tsmo}), $\sigma_{2s-1}([\delta_X,\wgts])$ is $s$ times a tensorial expression in $X, \dg, g$,
since $\partial_i\Gamma_{\nu j}^{\ell} = \frac{1}{3}(R_{i\nu j}^{\ \ \ \ell} +
R_{ij\nu }^{\ \ \ \ell})$ in normal coordinates.  The term with $\Gamma$ vanishes, so 
$\sigma_{2s-1}( [D_X,\wgts]) $ is $s$ times this tensorial expression.

The argument for $\sigma_{2s-1}([D_\cdot,\wgts]X$ is similar.  The term
with 
$\Gamma  $ vanishes.  
By (\ref{222}), (\ref{223a}), 
$$\sigma_{2s-1}([\delta_\cdot,\wgts]X)^a_b = 
i\sum_j\partial_t^j\partial_\xi^j|_{t=0, \xi=0} (p(\theta, \xi)^a_e\Gamma_{bc}^e(\gamma-t, \eta +\xi) X^c(\theta-t)).$$
By (\ref{tsmo}), the right hand side is linear in $s$ for Re$(s) <0$.  By (\ref{abc}), this implies 
the linearity in $s$ for Re$(s)>0.$  

Since $\sigma_{2s-1}([D_\cdot,\weight]X^*) = (\sigma_{2s-1}([D_\cdot, \weight]X))^*$, this 
symbol is also linear in $s$.
\end{proof}

\section{{\bf The Loop Group Case}}

In this section, we relate our work to Freed's work on based loop groups
$\Omega G$
\cite{Freed}.  We find a particular representation of the loop algebra that
controls 
the order of the curvature of the $H^1$ metric on $\Omega G.$

$\Omega G\subset LG$ has tangent space $T_\gamma\Omega G
= \{X\in T_\gamma LG: X(0) = X(2\pi) = 0\}$ in some Sobolev topology.  
Instead of using 
$D^2/d\gamma^2$ to define the Sobolev spaces, the usual choice is
$\Delta_{S^1} = -d^2/d\theta^2$ coupled to the
identity operator on the Lie algebra ${\mathfrak g}$.  Since this operator has
no kernel on $T_\gamma\Omega M$, 
$1 + \Delta$ is replaced by
 $\Delta$.  These changes in the $H^s$ inner product
do not alter the spaces of Sobolev sections, but the $H^s$ metrics on $\Omega G$ are 
no longer induced from a metric on $G$ as in the previous sections.

This simplifies the calculations of the Levi-Civita connections.
In particular,\\
 $[D_Z,\Delta^s] = 0$, so there is no term $A_XY$ as in (\ref{axy}).  
As a result, one can work directly with the six term formula (\ref{5one}).
For $X, Y, Z$ left
invariant vector fields, the first three terms on the right hand side of 
(\ref{5one}) vanish. Under the standing assumption that $G$ has a left
invariant, 
Ad-invariant inner product,
one obtains
$$2\nabla^{(s)}_XY = [X,Y] + \Delta^{-s}[X,\Delta^sY] +
\Delta^{-s}[Y,\Delta^sX]$$
\cite{Freed}.

It is an interesting question to compute the order of the curvature operator
as a function of $s$.  For based loops, Freed proved that this order is at
most $-1$.  In \cite{andres}, it is
shown that the order of $\Omega^s$ is at most $-2$ for all $s\neq 1/2, 1$ on
both $\Omega G$ and $LG$, and is exactly $-2$ for $G$ nonabelian.  
 For the case $s=1$, we have a much stronger result.

\begin{prop} The curvature of the
Levi-Civita connection for the $H^1$ inner product on $\Omega
G$ associated to $-\frac{d^2}{d\theta^2}\otimes {\rm Id}$ is a $\pdo$ of order $-\infty.$
\end{prop}

\noindent {\sc Proof:}
We give two quite different proofs.  

By \cite{Freed}, the $s=1$ curvature operator $\Omega = \Omega^{1}$
satisfies
$$\left\langle \Omega(X,Y)Z,W\right\rangle_1 = \left(\int_{S^1}[Y,\dot
Z],\int_{S^1}[X,\dot W]\right)_{\mathfrak g} - (X\leftrightarrow Y),$$
where  the inner product is the Ad-invariant form on the Lie algebra ${\mathfrak g}$.  We
want to write the right hand side of 
this equation as an $H^1$ inner product with $W$, in
order to recognize $\Omega(X,Y)$ as a $\pdo.$

Let $\{e_i\}$ be an orthonormal basis of ${\mathfrak g}$, considered
as a left-invariant frame of $TG$ and as global sections of $\gamma^*TG.$
 Let $\cc{i}{j}{k} = ([e_i,e_j],
e_k)_{\mathfrak g}$ be the structure constants of ${\mathfrak g}.$
(The Levi-Civita connection on left invariant vector fields
for the left invariant metric is
given by $\nabla_XY = \frac{1}{2}[X,Y]$, so the structure constants
are twice the Christoffel symbols.)  For $X = X^ie_i =
X^i(\theta)e_i, Y = Y^je_j,$ etc., integration by parts 
gives
$$\left\langle\Omega(X,Y)Z,W\right\rangle_1 = \left(\int_{S^1} \dot
Y^iZ^jd\theta\right)\left( \int_{S^1}\dot X^\ell W^m d\theta\right)
\cc{i}{j}{k}\cc{\ell}{m}{n}\delta_{kn} - (X\leftrightarrow Y).$$
Since
$$\int_{S^1}\cc{\ell}{m}{n}\dot X^\ell W^m =
\int_{S^1}\left(\delta^{mc}\cc{\ell}{c}{n}\dot X^\ell
e_m,W^be_b\right)_{\mathfrak g} = \left
\langle \Delta^{-1}(\delta^{mc}\cc{\ell}{c}{n} \dot X^\ell e_m),
W\right\rangle_1,$$
we get
\begin{eqnarray*}
\langle\Omega(X,Y)Z,W\rangle_1 &=& \left\langle
 \left[\int_{S^1} \dot Y^i Z^j\right]
\cc{i}{j}{k}\delta_{kn}\delta^{ms}\cc{\ell}{s}{n} \Delta^{-1}(\dot
X^\ell e_m),W\right\rangle_1- (X\leftrightarrow Y)\\
&=&\left\langle \left[ \int_{S^1}
a_j^k(\theta,\theta')Z^j(\theta')d\theta'\right] e_k,W\right\rangle_1,
\end{eqnarray*}
with
\begin{equation}\label{a}a_j^k(\theta,\theta') = \dot Y^i(\theta')
\cc{i}{j}{r}\delta_{rn}\delta^{ms}\cc{\ell}{s}{n} 
\left( \Delta
^{-1}( \dot X^\ell
e_m)\right)^k(\theta) - (X\leftrightarrow Y).
\end{equation}

We now show that $Z\mapsto \left(\int_{S^1}
a_j^k(\theta,\theta')Z^j(\theta')d\theta'\right)e_k$ is a smoothing
operator.  Applying Fourier transform and Fourier inversion to $Z^j$
yields
\begin{eqnarray*} \int_{S^1} a_j^k(\theta,\theta')Z^j(\theta')d\theta'
&=& \int_{S^1\times\R\times S^1}
a_j^k(\theta,\theta')e^{i(\theta'
-\theta'')\cdot\xi}Z^j(\theta'')d\theta''d\xi d\theta'\\
&=&
\int_{S^1\times\R\times S^1} \left[ a_j^k(\theta,\theta')e^{-i(\theta
-\theta')\cdot\xi}\right]e^{i(\theta
-\theta'')\cdot\xi}Z^j(\theta'')d\theta''d\xi d\theta',
\end{eqnarray*}
so $\Omega(X,Y)$ is a $\pdo $ with symbol 
\begin{equation}\label{b} b_j^k(\theta,\xi) =
\int_{S^1} a_j^k(\theta,\theta') \eee d\theta',
\end{equation}
with the usual mixing of local and global notation.

For fixed $\theta$,
(\ref{b})  contains the Fourier transform of $\dot Y^i(\theta')$  and $\dot X^i(\theta')$, as
these are the only $\theta'$-dependent terms in (\ref{a}).
Since the
Fourier transform is taken in a local chart with respect to a
partition of unity, and since in each chart $\dot Y^i$ and $\dot X^i$ times the
partition of unity function is compactly supported, the Fourier
transform of $a_j^k$ in each chart is rapidly decreasing.  Thus
$b_j^k(\theta,\xi)$ is the product of a rapidly decreasing function
with $e^{i\theta\cdot\xi}$, and hence is of order $-\infty.$

We now give a second proof.  For all $s$,
$$\nabla_X Y = \frac{1}{2}[X,Y] -\frac{1}{2} \Delta^{-s}[\Delta^sX,Y]
+\frac{1}{2}\Delta^{-s}[X,\Delta^sY].$$
Label the terms on the right hand side (1) -- (3).
 As an operator on $Y$ for fixed $X$, the symbol of (1) is
$\sigma((1))^a_\mu = \frac{1}{2}X^ec_{\e\mu}^a.$
 Abbreviating $\xii{-s}$ by $\xi^{-2s}$, we have
\begin{eqnarray*} \sigma((2))^a_\mu &\sim & -\frac{1}{2}c_{\e\mu}^a
\left[ \xi^{-2s}\Delta^sX^\e -\frac{2s}{i}\xi^{-2s-1}
\partial_\theta\Delta^s X^\e  \right.\\
&&\ \ \  \left. +\sum_{\ell=2}^\infty\frac{(-2s)(-2s-1)
\ldots(-2s-\ell+1)}{i^\ell \ell!}\xi^{-2s-\ell}
\partial_\theta^\ell\Delta^s X^\e \right]\\
\sigma((3))^a_\mu &\sim &  \frac{1}{2}c_{\e\mu}^a
\left[ X^\e+ \sum_{\ell=1}^\infty \frac{(-2s)(-2s-1)
\ldots(-2s-\ell+1)}{i^\ell \ell!} \xi^{-\ell}\partial_\theta^\ell X^\e\right].
\end{eqnarray*}
Thus
\begin{eqnarray}\label{fourone}
\sigma(\nabla_X)^a_\mu  &\sim& \frac{1}{2}c_{\e\mu}^a\left[ 2X^\e
   -\xi^{-2s}\Delta^sX^\e
+\frac{2s}{i}
\xi^{-2s-1}\partial_\theta\Delta^sX^\e\right. \nonumber\\
&&\ \ \
   -\sum_{ \ell=2}^\infty\frac{(-2s)(-2s-1)\ldots(-2s-\ell+1)}{i^\ell \ell!}
\xi^{-2s-\ell}\partial_\theta^\ell\Delta^s X^\e \\
&&\ \ \  \left. + \sum_{\ell=1}^\infty \frac{(-2s)(-2s-1)
\ldots(-2s-\ell+1)}{i^\ell \ell!} \xi^{-\ell}\partial_\theta^\ell
 X^\e. \right].\nonumber
\end{eqnarray}

Set $s=1$ in (\ref{fourone}), and replace $\ell$
by
$\ell-2$ in the first infinite sum.  Since $\Delta = -\partial_\theta^2$, a
little algebra gives
\begin{equation}\label{fourtwo}
\sigma(\nabla_X)^a_\mu \sim c_{\e\mu}^a\sum_{\ell=0}^\infty
\frac{(-1)^\ell}{i^\ell}
\partial_\theta^\ell X^\e\xi^{-\ell}
=  \ad\left( \sum_{\ell=0}^\infty
\frac{(-1)^\ell}{i^\ell}\partial_\theta^\ell
X\xi^{-\ell}
\right).
\end{equation}

Denote the infinite sum in the last term of (\ref{fourtwo})
by $W(X,\theta,\xi)$. The map
$X\mapsto W(X,\theta,\xi)$ takes the  Lie algebra of left invariant vector
fields on $LG$ to the Lie algebra
$L{\mathfrak g}[[\xi^{-1}]], $
the space of formal $\pdo$s of nonpositive integer order on the trivial bundle
$S^1\times{\mathfrak g} \to S^1$, where the Lie bracket on the
target involves multiplication of power series and bracketing in
${\mathfrak g}.$  We claim that this map is a Lie algebra homomorphism.
Assuming this, we see that
\begin{eqnarray*} \sigma\left(\Omega(X,Y)\right) &=&
  \sigma\left([\nabla_X,\nabla_Y] -\nabla_{[X,Y]}\right)
\sim \sigma\left( [\ad W(X), \ad W(Y)] - \ad W([X,Y]) \right)\\
&=& \sigma\left( \ad ( [W(X), W(Y)]) - \ad W([X,Y]) \right) = 0,
\end{eqnarray*}
which proves that $\Omega(X,Y)$ is a smoothing operator.

To prove the claim,
set $X = x^a_n\eff{n}e_a, Y =y^b_m\eff{m}e_b$. 
Then
\begin{eqnarray*} W([X,Y]) &=&
 W( x^ny^m\eff{(n+m)}c_{ab}^k e_k) =\sum_{\ell=0}^\infty \frac{(-1)^\ell}
{i^\ell } c_{ab}^k \partial_\theta^\ell
 \left(x^a_ny^b_m\eff{(n+m)}\right) \xi^{-\ell}e_k\\
  {[} W(X)  ,   W(Y)]
&=& \sum_{\ell=0}^\infty \sum_{p+q = \ell}
\frac{(-1)^{p+q}}{i^{p+q}} \partial_\theta^p \left(
x^a_n\eff{n}\right) \partial_\theta^q
\left( y^b_m\eff{m}\right)\xi^{-(p+q)}c_{ab}^k e_k,
\end{eqnarray*}
and these two sums are clearly equal.
\hfill $\Box$

\bigskip

It would be interesting to understand how the map $W$ fits into the
representation theory of the loop algebra $L{\calg}.$ 
\bigskip

\large
\noindent {{\bf Part II. Characteristic Classes on $LM$}}
\normalsize
\bigskip

In this part, we construct a general theory of Chern-Simons
classes on certain infinite rank bundles including the frame/tangent bundle of 
loop spaces,
following the construction of primary characteristic classes
in \cite{P-R2}. The primary  classes vanish on the tangent bundles of
loop spaces, which forces the
consideration of secondary classes. 
The key ingredient is to replace the ordinary matrix trace in the Chern-Weil
theory of
 invariant polynomials
on  finite dimensional Lie groups with the Wodzicki residue on invertible bounded
$\pdo$s.  

As discussed in the Introduction, there are absolute and relative versions of Chern-Simons theory.  We use the relative version, which assigns an odd degree form to a pair 
of connections.
In particular, for $TLM$, we can use the $L^2$ (i.e. $s=0$) and
 $s=1$ Levi-Civita connections to form Wodzicki-Chern-Simons (WCS) classes associated to a metric on $M$. 
 
 In \S\ref{CSCLS}, we develop the general theory of Wodzicki-Chern and WCS classes for 
 bundles with structure group $\pdo_0^*$, the group of invertible classical zeroth order pseudodifferential operators.   We show the vanishing of 
  the Wodzicki-Chern classes of $LM$ and more general mapping spaces. 
As in finite dimensions, we show the existence of  WCS classes in 
  $H^n(LM,\BbC)$ if dim$(M) = n$ is odd (Definition \ref{def:WCS})
  and give the local expression for the WCS classes associated to the Chern character
  (Theorem \ref{thm:5.5}).
   In Theorem \ref{WCSvan}, we prove that the
    Chern character WCS class vanishes if dim$(M) \equiv 3
   \ ({\rm mod}\ 4)$.
In \S\ref{dimfive}, we associate to every circle action $a:S^1\times M^n\to M^n$
 an $n$-cycle $[a]$
 in $LM$.  For a specific metric on $S^2\times S^3$ and a specific circle action $a,$
 we prove via exact computer calculations that the WCS class is nonzero by integrating it over $[a].$
 Since the corresponding integral for the cycle associated to the trivial action 
 is zero, $a$ cannot be homotoped to the trivial action. 
We use this result to prove that $\pi_1({\rm Diff}
 (S^2\times S^3))$ is infinite.

Throughout this part, $H^*$ always refers to de Rham cohomology for complex valued forms.  By \cite{beggs}, $H^*(LM)\simeq H^*_{\rm sing}(LM,\BbC).$

\section{{\bf Chern-Simons Classes on Loop Spaces}}\label{CSCLS}

We begin in \S5.1 with a review of Chern-Weil and Chern-Simons theory in
finite dimensions, following  \cite{C-S}.  

In
\S5.2, we discuss Chern-Weil and Chern-Simons theory on a class of infinite rank bundles
including the frame bundles of loop spaces.  As in \S2.7,  the geometric
structure group of these bundles
 is $\pdo_0^*$, so we need a trace on the Lie algebra
$\pdo_{\leq 0}$ to define invariant polynomials.  There are two
types of  traces, one given by taking the zeroth order symbol and one given by
the Wodzicki residue \cite{paycha-lescure}, \cite{ponge}.  Here we only consider the 
Wodzicki residue trace.

Using this trace, we generalize the usual definitions of Chern and Chern-Simons classes in
de Rham cohomology.  In particular,
given a $U(n)$-invariant polynomial $P$ of degree $k$, we define a corresponding 
WCS class $CS^W_P\in H^{2k-1}(LM)$ if dim$(M) = 2k-1.$    We are forced to consider these secondary classes, because the Wodzicki-Chern classes of mapping spaces
$\maps(N,M)$ vanish.  In Theorem \ref{thm:5.5}, we give an exact expression for the WCS
classes associated to the Chern character.
In Theorem \ref{WCSvan}, we show that these WCS classes in $H^{4k+3}(LM^{4k+3})$
vanish; in contrast, in finite dimensions, the Chern-Simons classes associated to the Chern character vanish in $H^{4k+1}(M^{4k+1}).$

\subsection{{\bf Chern-Weil  and Chern-Simons Theory for Finite Dimensional
  Bundles} }

We first review the Chern-Weil construction. 
Let $G$ be a finite dimensional Lie group with Lie algebra $\g$, and let
 $G\to F\to M$ be a principal $G$-bundle over a manifold $M$. 
Set $
 \g^k=\g^{\otimes k}$ and let
\begin{equation*}I^k(G)
= \{P:\g^k\to \C\ | P\ \text{symmetric,
  multilinear, Ad-invariant}\}
\end{equation*}
be the degree $k$ Ad-invariant polynomials on $\g.$

\begin{rem}
For classical Lie groups $G$, $I^k(G)$ is generated by the polarization of
the Newton polynomials $\Tr(A^\ell)$, where $\Tr$ is the usual trace on finite
dimensional matrices.
\end{rem}

For $\phi\in\Lambda^\ell(F,\g^k)$, $P\in I^k(G)$, set
$P(\phi)=P\circ \phi\in\Lambda^\ell(F)$.

\begin{thm}[The Chern-Weil Homomorphism \cite{K-N}] \label{previous}
Let $F\to M$ have a connection $\theta$ with curvature $\Omega_F\in
\Lambda^2(F,\g)$. For $P\in I^k(G)$, $P(\Omega_F)$ is a closed
 invariant real form on $F$, and so
determines a closed form
$P(\Omega_M)\in \Lambda^{2k}(M)$.
The Chern-Weil map
\begin{equation*}
\oplus_{k}I^k(G)\to H^{*}(M), \ P\mapsto [P(\Omega_M)]
\end{equation*}
is a well-defined algebra homomorphism, and in particular is independent of the choice of
connection on $F$.
\end{thm} 

The proof depends on:
\begin{itemize}
\item (The {\it commutativity property}) 
For $\phi\in\Lambda^{\ell}(F,\g^k)$, 
\begin{equation}\label{eq:deri}
d(P(\phi))=P(d\phi).
\end{equation}
\item (The {\it infinitesimal invariance property})
For $\psi_i\in\Lambda^{\ell_i}(F,\g)$, $\phi\in\Lambda^{1}(F,\g)$ and $P\in
  I^k(G)$, 
\begin{equation}\label{eq:inva}
\sum^k_{i=1}
(-1)^{\ell_1+\dots+\ell_i}P(\psi_1\wedge\dots\wedge[\psi_i,\phi]\wedge\dots
\psi_l)=0. 
\end{equation}
\end{itemize}
$[P(\Omega_M)]$ is
called the {\it characteristic class} of $P$. For example, the characteristic class
 associated to  $\Tr(A^k)$ is the k${}^{\rm th}$ component of the Chern character of $F$.

Part of the theorem's content
is that for any two connections on $F$,
$P(\Omega_1) - P(\Omega_0) = 
dCS_P(\theta_1,\theta_0)$  
for some odd form $CS_P(\nabla_1, \nabla_0)$.  Explicitly, 
\begin{equation}\label{5.1}
CS_P(\theta_1,\theta_0) = \int_0^1 P(\theta_1-\theta_0,\overbrace{\Omega_t,...,\Omega_t}^{k-1})
\ dt
\end{equation}
where 
$$\theta_t = t\theta_0+(1-t)\theta_1,\ \ \Omega_t = d\theta_t+\theta_t\wedge\theta_t$$ \cite[Appendix]{chern}.  

\begin{rem}
For $F\stackrel{\pi}{\to} M$,  
$\pi^*F\to F$ is trivial.
Take $\theta_1$ to be the flat connection on $\pi^*F$
with respect to a fixed trivialization.
Let $\theta_1$ also
denote the connection $\chi^*\theta_1$ on $F$, 
where $\chi$ is the global section of $\pi^*F.$  For any other connection $\theta_0$ on $F$,  $\theta_t = t\theta_0, \Omega_t = t\Omega_0 + (t^2-t)\theta_0\wedge \theta_0$. 
 Assume an invariant  polynomial $P$ takes values in $\R.$  Then  we
obtain the formulas for the transgression form $TP(\Omega_1)$
on $F$: for 
\begin{equation}\label{eq:ChernSimons}
\phi_t =t\Omega_1+\frac{1}{2}(t^2-t)[\theta,\theta],\ \ 
TP(\theta)=l\int_0^1 P(\theta\wedge \phi^{k-1}_t)dt,
\end{equation}
$dTP(\theta)=P(\Omega_1)\in \Lambda^{2l}(F)$
\cite{C-S}.  $TP(\Omega_1)$ pushes down to an $\BbR/\BbZ$-class on $M$,
the absolute Chern-Simons class.
\end{rem}

As usual, these formulas carry over to connections $\nabla = d+\omega$
on vector bundles $E\to M$ in the form
\begin{equation}\label{5.11}
CS_{P}(\nabla_1,\nabla_0) = \int_0^1 P(\omega_1-\omega_0,\Omega_t,..., 
\Omega_t)\ dt,
\end{equation}
since $\omega_1-\omega_0$ and
$\Omega_t$ are globally defined forms.

 \subsection{{\bf Chern-Weil and Chern-Simons Theory for $\pdo_0^*$-Bundles}}

 Let $\mathcal E\to\mathcal M$ be an infinite rank bundle
 over a paracompact Banach manifold
 $\mathcal M$, with the fiber of $\mathcal E$ 
 modeled on a fixed Sobolev class of sections of 
 a finite rank hermitian vector  bundle $E\to N$, and with structure group $\GG(E)$.  
 For such $\GG$-bundles,
 we can produce
 primary and secondary characteristic classes 
 once we choose a trace on $\pdo_{\leq 0}(E)$.
 Since the adjoint action of $\GG$ on $\pdo_{\leq 0}$ is by conjugation,  a trace on $\pdo_{\leq 0}$ will extend to a polynomial on forms
satisfying (\ref{eq:deri}), (\ref{eq:inva}), so the finite dimensional proofs extend.  
 
These traces were classified  in \cite{lesch-neira, paycha-lescure}, although there are slight variants
in our special case $N= S^1$ \cite{ponge}.  Roughly speaking, the traces fall into two classes, the leading order symbol trace \cite{P-R2} and the Wodzicki residue.  In this paper,
we consider only the Wodzicki residue, and refer to \cite{lrst} for the leading order symbol
trace.

For simplicity, we mainly restrict to the generating invariant polynomials $P_k(A) = A^k$, and 
only consider $\mathcal E = TLM$, which we recall is the complexified tangent bundle.  We will work with vector bundles rather than principal bundles.

\begin{defn}  \label{def:WCS} 
(i) The k${}^{\rm th}$ {\it Wodzicki-Chern (WC) form} of a $\pdo_0^*$-connection
$\nabla$ on $TLM$ with curvature $\Omega$ is
\begin{equation}\label{5.1a}
c_k^W(\Omega)(\gamma) =\frac{1}{k!}
 \int_{S^*S^1}\tr\sigma_{-1}(\Omega^{k}) \ d\xi  dx.
\end{equation}
Here we recall that for each $\gamma\in LM$,
$\sigma_{-1}(\Omega^k)$ is a $2k$-form with values in  endomorphisms
 of a trivial bundle
over $S^*S^1$.

(ii) The k${}^{\rm th}$ {\it Wodzicki-Chern-Simons (WCS) form} of two $\pdo_0^*$-connections 
$\nabla_0,\nabla_1$ on $TLM$ is
\begin{eqnarray}\label{5.22}
CS^W_{2k-1}(\nabla_1,\nabla_0) &=&\frac{1}{k!}
 \int_0^1 \int_{S^*S^1}\tr\sigma_{-1}((\omega_1-\omega_0)\wedge 
(\Omega_t)^{k-1})\ dt\\ 
&=&\frac{1}{k!} \int_0^1 {\rm res}^{\rm w} 
[(\omega_1-\omega_0)\wedge 
(\Omega_t)^{k-1}]\ dt.\nonumber
\end{eqnarray}

(iii) The  k${}^{\rm th}$ {\it Wodzicki-Chern-Simons form} associated to a Riemannian metric 
$g$ 
on $M$, denoted $CS^W_{2k-1}(g)$,  is $CS^W_{2k-1}(\nabla_1,\nabla_0)$, where $\nabla_0, \nabla_1$ refer to the 
$L^2$ and $s=1$ Levi-Civita connections on $LM$, respectively.

(iv) Let $\Sigma = \{\sigma\}$ be the group of permutations of $\{1,...,k\}$. Let $I:
1\leq i_1< ...< i_\ell = k$ be a partition of $k$ (i.e. with $i_0=0$, $\sum_{j=1}^k
 (i_j-i_{j-1}) = k$) .  For  the symmetric, $U(n)$-invariant, 
multilinear form on ${\mathfrak u}(n)$
\begin{eqnarray*} P_I(A_1,A_2,...,A_k) &=& \frac{1}{k!}
\sum_\sigma \tr(A_{\sigma(1)}\cdot...\cdot A_{\sigma(i_1)})
\tr(A_{\sigma(i_1+1)}\cdot...\cdot A_{\sigma(i_2)})\\
&&\qquad \cdot ...\cdot \tr(A_{\sigma(i_{\ell-1})}
\cdot ...\cdot A_{\sigma(k)}),
\end{eqnarray*}
define the symmetric, $\GG$-invariant, multilinear form on $\pdo_{\leq 0}$ by
\begin{eqnarray*} P_I^W(B_1,...,B_k) &=&  \frac{1}{k!}
 \sum_\sigma\left( \intss \tr\sigma_{-1}(B_{\sigma(1)}\cdot...\cdot B_{\sigma(i_1)}) \right.\\
&&\qquad \left. \cdot
 \intss \tr\sigma_{-1}
(B_{\sigma(i_1+1)}\cdot...\cdot B_{\sigma(i_2)})\right. \\
&&\qquad \left. \cdot ...\cdot \intss \tr\sigma_{-1}(B_{\sigma(i_{\ell-1})} 
\cdot ...\cdot B_{\sigma(k)})\right).
\end{eqnarray*}
  The {\it Wodzicki-Chern form associated
to $P_I$} for a $\pdo_0^*$-connection on $TLM$ with curvature $\Omega$ is 
\begin{eqnarray}\label{wcpi} c_{P_I}^W(\Omega) &=& 
P_I^W(\Omega,\Omega,...,\Omega)\\
 &=& \frac{1}{k!}
\intss \tr\sigma_{-1}(\Omega^{k_1}) \cdot \intss \tr\sigma_{-1}(\Omega^{k_2})
\cdot
...\cdot \intss \tr\sigma_{-1}(\Omega^{k_\ell}\nonumber )\\
&=& \frac{k_1!k_2!\cdot...\cdot k_\ell !}{k!} 
c_{k_1}^W(\Omega)c_{k_2}^W(\Omega)\cdot ...\cdot c_{k_\ell}^W(\Omega),\nonumber
\end{eqnarray}
where $k_1=i_1-i_0, k_2 = i_2-i_1,...,
k_\ell = i_\ell - i_{\ell-1}$.

Setting $K = (k_1,...,k_\ell)$, we also denote $c_{P_I}^W(\Omega)$ by $c_K^W(\Omega).$


(v) 
Let $\nabla_0,\nabla_1$ be $\pdo_0^*$-connections on $TLM$ with connection forms
$\omega_0, \omega_1,$ respectively.  The 
 {\it Wodzicki-Chern-Simon form associated to $P_I$ and $\nabla_0, \nabla_1$}
is
$$
CS^W_{P_I}(\nabla_1,\nabla_0) = \int_0^1 P_I^W(\omega_1-\omega_0,\Omega_t,...,
\Omega_t)dt. $$ 
\end{defn}

In (iv) and (v), we do not bother with a normalizing constant, since we do not claim that 
there is a normalization which gives classes with integral periods.  
Note that the k${}^{\rm th}$ WCS class is associated to $P_k(A_1,...,A_k) = \tr(A_1\cdot
...\cdot A_k)$, i.e. the partition $K = (k)$, or in other words to the polynomial giving the 
k${}^{\rm th}$ component of the Chern character.

As in finite dimensions, $c_k^W(\nabla)$ is a closed $2k$-form, with de Rham cohomology
class $c_k(LM)$
 independent of $\nabla$, as $c_k^W(\Omega_1)  - c_k^W(\Omega_0) =
dCS^W_{2k-1}(\nabla_1,\nabla_0).$  
 
\begin{rem}  It is an interesting question to determine all the $\pdo_0^*$-invariant polynomials on 
$\pdo_{\leq 0}.$  As above, $U(n)$-invariant polynomials combine with the Wodzicki residue 
(or the other traces on $\pdo_{\leq 0}$) to give $\pdo_0^*$-polynomials,
but there may be others.  
\end{rem}

The tangent space $TLM$, and more generally mapping spaces
Maps$(N,M)$  with $N$ closed
have vanishing Wodzicki-Chern classes.  Here we take a Sobolev topology on Maps$(N,M)$ for some
large Sobolev parameter, so that Maps$(N,M)$ is a paracompact Banach manifold.
We denote the de Rham class of $c_{P_I}^W(\Omega)$ for a connection on $\mathcal E$ by
$c_{P_I}(\mathcal E).$

\begin{prop} \label{prop:maps} Let $N, M$ be closed manifolds, and let {\rm Maps}${}_f(N,M)$ denote
the component of a fixed $f:N\to M$.  Then the  cohomology classes
$c_{P_I}^W({\rm Maps}_f(N,M)) $ of $T{\rm Maps}(M,N)$ vanish.
\end{prop}

\begin{proof}
For $TLM$, the $L^2$ connection in Lemma \ref{lem:l2lc}
has curvature $\Omega$ which is a multiplication operator.  Thus $\sigma_{-1}(\Omega)$ and hence $\sigma_{-1}(\Omega^{i})$ are zero, 
so the WC forms $c_{P_I}(\Omega)$ also vanish.

For $n\in N$ and $h:N\to M$,
let $\ev_n: {\rm Maps}_f(N,M)$ be $\ev_n(h) = h(n).$ 
 Then $D_XY(h)(n) \stackrel{\rm def}{=} 
 (\ev_h^*\nabla^{LC,M})_XY(h)(n)$ is the $L^2$ Levi-Civita connection on \\
Maps$(N,M).$
As in Lemma \ref{lem:l2lc},
the curvature of $D$ is a 
a multiplication operator.  Details are left to the reader.
\end{proof}

\begin{rem}  (i) These mapping spaces fit into the framework of the Families Index Theorem in 
the case of a trivial fibration
$Z\to M\stackrel{\pi}{\to} B$ of closed manifolds.  Given a
finite rank bundle $E\to M$, we get an associated infinite rank bundle $\calE = \pi_*E
\to B$.  For the fibration $N\to N\times {\rm Maps}(N,M)\to {\rm Maps}(N,M)$ and $E = {\rm
  ev}^*TM$, $\calE$ is $T{\rm Maps}(N,M).$  A connection $\nabla$ on $E$ induces a connection $\nabla^{\calE}$ on $\calE$ defined by
\begin{equation*}
(\nabla^{\calE}_Z s)(b)
(z)=\left( (\ev^*\theta^u)_{(Z,0)} u_s\right)(b,z).
\end{equation*}
Here $u_s(b,z)=s(b)(z)$.
The curvature $\Omega^{\calE}$ satisfies
\begin{equation*}\label{eq:pullback}
\Omega^{\calE}(Z,W)s(b)(z)=(\ev^*\Omega ) ((Z,0),(W,0)) u_s(b,z).
\end{equation*}
This follows from
\begin{equation*}
\Omega^{\calE}(Z,W)s(b)(z)= [\nabla^{\calE}_Z \nabla^{\calE}_W
-\nabla^{\calE}_W \nabla^{\calE}_Z -\nabla^{\calE}_{[Z,W]}]
s(b)(z).
\end{equation*}
Thus the connection and curvature forms take values in multiplication operators, and 
so $c_k^W(\calE) = 0.$

If the fibration is nontrivial, the connection on $\calE$ depends on the choice of a horizontal complement to $TZ$ in $TM$, and the corresponding connection and curvature forms take
values in first order differential operators.  

(ii) In finite dimensions, odd Chern forms of complexified real bundles like\\
 $T{\rm Maps}(N,M)$ vanish, because the form involves a composition of an odd number of skew-symmetric matrices.  In contrast, odd WC forms involve terms like
$\sigma_{-1}(\Omega^1)\wedge\Omega^M\wedge...\wedge\Omega^M,$ where $\Omega^1$ is the curvature of the $s=1$ Levi-Civita connection.  By Lemma 
\ref{old2.2}(ii), $\sigma_{-1}(\Omega^1)$ is not skew-symmetric as an endomorphism.  Thus
it is not obvious that the odd WC forms vanish.

Similarly, in finite dimensions the Chern-Simons form for the odd Chern classes of complexified real bundles vanish, but this need not be the case for WCS forms.  In fact, we will produce nonvanishing 
WCS classes associated to $c_3^W(TLM^5)$ in \S\ref{dimfive}.

\end{rem}

In finite dimensions, Chern classes are topological obstructions to the
reduction of the structure group and geometric obstructions to the existence
of a flat connection.  
Wodzicki-Chern classes for $\pdo_0^*$-bundles 
are also topological and geometric obstructions, but
the geometric information is a little more refined due to the grading on the
Lie algebra 
 $\pdo_{\leq 0}$.

\begin{prop}
 Let $\calE\to\calB$ be an infinite rank $\GG$-bundle, for
  $\GG$ acting on
$E\to N^n$.  
If $\calE$ admits a reduction to the gauge group $\calG(E)$, then
  $c_k^W(\calE)  =   0$ for all $k$, and hence $c_{P_I}^W(\calE) =0$ for all $P_I$.
If $\calE$ admits a 
   $\GG$-connection whose
  curvature has order $-k$, then $
  c_{\ell}(\calE) =0$ for $\ell \geq [n/k].$
 \end{prop}

\begin{proof}  If the structure group of  $\calE$ reduces to the gauge
  group, there exists a connection one-form
  with values in Lie$(\calG) = {\rm End}(E)$, the Lie algebra of multiplication
  operators.  Thus the Wodzicki residue of powers of the curvature vanishes,
  so the Wodzicki-Chern classes vanish.
For the second statement, the order of the curvature is less than
$-n$ for $\ell \geq [n/k]$,  so the Wodzicki residue
  vanishes in this range. 
  \end{proof}
  
  However, we do not have examples of nontrivial WC classes; cf.~\cite{lrst}, where it is 
  conjectured that these classes always vanish.  
  \bigskip

The relative WCS form is not difficult to compute.  

\begin{prop}  Let $\sigma$ be in the group of permutations of $\{1,\ldots,2k-1\}.$ Then
\begin{eqnarray}\label{5.4}
\lefteqn{CS^W_{2k-1}(g)(X_1,...,X_{2k-1}) }\\
&=&
\frac{2}{(2k-1)!} \sum_{\sigma} {\rm sgn}(\sigma) \int_{S^1}\tr [
(-2R(X_{\sigma(1)},\dg)
-R(\cdot,\dg)X_{\sigma(1)} + R(X_{\sigma(1)},\cdot)\dg) \nonumber\\
&&\qquad 
\cdot (\Omega^M)^k(X_{\sigma(2)},..X_{\sigma(2k-1)} )].\nonumber
\end{eqnarray}
\end{prop}

\begin{proof}
$$\sigma_0((\omega_1-\omega_0)_X)^a_b = \Gamma_{cb}^aX^c -\Gamma_{cb}^aX^c = 0.$$
Thus 
\begin{equation}\label{cswint}
CS^W_{2k-1}(g) = \int_0^1 \int_{S^*S^1}\tr\sigma_{-1}(\omega_1-\omega_0)\wedge (\sigma_0(\Omega_t))^k\ dt.
\end{equation}
Moreover,
\begin{eqnarray*}\sigma_0(\Omega_t) &=& td(\sigma_0(\omega_0)) + (1-t)d(\sigma_0(\omega_1)) \\
&&\qquad + 
(t\sigma_0(\omega_0) + (1-t)\sigma_0(\omega_1))\wedge (t\sigma_0(\omega_0) + (1-t)\sigma_0(\omega_1))\\
&=& d\omega^M + \omega^M\wedge \omega^M\\
&=& \Omega^M.
\end{eqnarray*}
Therefore
\begin{equation}\label{5.3}
CS^W_{2k-1}(g) = \int_0^1 \int_{S^*S^1}\tr [\sigma_{-1}(\omega_1)
\wedge (\Omega^M)^k]\ dt,
\end{equation}
since $\sigma_{-1}(\omega_0) = 0.$  We can drop the integral over $t$. 
The integral over the $\xi$ variable contributes a factor of $2$: the integrand has
a factor of $|\xi|^{-2}\xi$, which equals $\pm 1$ on the two components of $S^*S^1$.
Since the fiber of $S^*S^1$ at a fixed $\theta$ consists of two points 
with opposite orientation, the ``integral" over each fiber is $1-(-1) = 2.$ 
Thus
\begin{eqnarray}\label{5.4a}  \lefteqn{
CS^W_{2k-1}(g)(X_1,...X_{2k-1})    }  \\
&=& =  \frac{2}{(2k-1)!}   \sum_\sigma {\rm sgn}(\sigma)  \int_{S^1}\tr[
(-2R(X_{\sigma(1)},\dg)
-R(\cdot,\dg)X_{\sigma(1)} + R(X_{\sigma(1)},\cdot)\dg)\nonumber\\
&&\qquad
\cdot (\Omega^M)^k(X_{\sigma(2)},..X_{\sigma(2k-1)} )]\nonumber
\end{eqnarray}
by Lemma \ref{old2.1}.
\end{proof}

This produces odd classes in the de Rham cohomology of the loop space of an odd
dimensional manifold.

\begin{thm}\label{thm:5.5}
 (i)  Let dim$(M) = 2k-1$ and let $P$ be a $U(n)$-invariant polynomial of degree 
$k.$  Then $c^W_P(\Omega) \equiv 0$ for any $\pdo_0^*$-connection $\nabla$ on
 $TLM.$  Thus $CS^W_P(\nabla_1,\nabla_0)$ is closed and defines a 
 class $[CS^W_P(\nabla_1,\nabla_0)]\in H^{2k-1}(LM).$  In particular, we can
 define $[CS^W_P(g)]\in H^{2k-1}(LM)$ for a Riemannian metric $g$ on $M$.

(ii)  For dim$(M) = 2k-1$, the k${}^{\it th}$ Wodzicki-Chern-Simons form $CS^W_{2k-1}(g)$
simplifies to 
 \begin{eqnarray}\label{csg}
\lefteqn{CS^W_{2k-1}(g)(X_1,...,X_{2k-1}) }\nonumber \\
&=&
\frac{2}{(2k-1)!} \sum_{\sigma} {\rm sgn}(\sigma) \int_{S^1}\tr[
(-R(\cdot,\dg)X_{\sigma(1)} + R(X_{\sigma(1)},\cdot)\dg)\\
&&\qquad 
\cdot (\Omega^M)^{k-1}(X_{\sigma(2)},..X_{\sigma(2k-1)} )].\nonumber
\end{eqnarray}

 \end{thm}
 
 \begin{proof} (i) Let $\Omega$ be the curvature of $\nabla.$
 $c^W_P(\Omega)(X_1,\dots, X_{2k})(\gamma)$ is a sum of monomials of the form
(\ref{wcpi}).  This is
a $2k$-form on $M$, and hence vanishes.  

 (ii) Since
 $$R(X_{1},\dg)
\cdot (\Omega^M)^k(X_{2},..X_{2k-1}) = 
[i_{\dg}\tr(\Omega^{k})](X_1,...X_{2k-1}) = \tr(\Omega^k)(\dg, X_1,\ldots,X_{2k-1}),$$
the first term on the right hand side of (\ref{5.4a}) vanishes on a $(2k-1)$-manifold.

\end{proof}

\begin{rem} There are several variants to the construction of relative WCS classes.

(i) If we define the transgression form $Tc_k(\nabla)$ with the Wodzicki residue
replacing the trace in (\ref{eq:ChernSimons}), it is easy to check that $Tc_k(\nabla)$
involves $\sigma_{-1}(\Omega).$  For $\nabla$ the $L^2$ connection, this WCS class vanishes.  For $\nabla$ the $H^s$ connection, $s>0$, $\sigma_{-1}(\Omega)$ involves 
the covariant derivative of the curvature of $M$ (cf.~Lemma \ref{old2.2} for $s=1.$)  Thus the
relative WCS class is easier for computations than the absolute class $[Tc_k(\nabla)].$

(ii) If we define $CS_k^W(g)$ using the Levi-Civita connection for the $H^s$ 
metric instead of
the $H^1$ metric, the WCS class is simply multiplied by the artificial parameter $s$ by 
Lemma \ref{lem:33}.  Therefore setting $s=1$ is not only computationally convenient, it 
regularizes the WCS, in that it extracts the $s$-independent information.
This justifies the following definition:

\begin{defn}  \label{def:regularized}
The {\it regularized  $k^{th}$ WCS class} associated to a Riemannian metric 
$g$ on $M$ is $CS_k^{W, {\rm reg}}(g) \stackrel{\rm def}{=} 
CS_k^W(\nabla^{1},\nabla^0)$, where $\nabla^{1}$ is the $H^1$ connection
and $\nabla^0$ is the $L^2$ Levi-Civita connection.  
\end{defn}  

\end{rem}

\bigskip

We conclude this section with a vanishing result that does not have a finite dimensional 
analogue.
  \begin{thm} \label{WCSvan}
  The {\it  k}${}^{\it  th}$ WCS class $CS_k^W(g)$
   vanishes if ${\rm dim}(M) \equiv 3 \ ({\rm mod}\  4).$
  \end{thm}

\begin{proof}  Let dim$(M) =2 k-1$. Since $\Omega^M$ takes values in skew-symmetric endomorphisms, 
so does
$(\Omega^M)^{k-1}$ if $k$ is even, i.e. if ${\rm dim}(M) \equiv 3 \ ({\rm mod}\  4).$
  The term 
  $-R(\cdot,\dg)X_{\sigma(1)} + R(X_{\sigma(1)},\cdot)\dg$ in (\ref{csg}) is a symmetric
  endomorphism.  For in Riemannian normal coordinates, this term is
  $(-R_{bdca} +R_{cbda})X^c\dg^d \equiv A_{ab}$, say, so the curvature terms in 
  $A_{ab} - A_{ba}$ are
  \begin{eqnarray*}
  -R_{bdca} +R_{cbda} + R_{adcb} - R_{cadb} &=& -R_{bdca} +R_{cbda} 
  + R_{cbad} - R_{dbca}\\
  &=& -R_{bdca} +R_{cbda} -R_{cbda} + R_{bdca}=0.
  \end{eqnarray*}
Thus the integrand in  (\ref{csg}) is the trace of a symmetric endomorphism composed with a skew-symmetric endormorphism, and so 
vanishes.
\end{proof}

\begin{exm}  We contrast Theorem \ref{WCSvan} with the situation in finite dimensions. 
Let dim$(M)=3.$ 
The only invariant monomials of degree two are $\tr(A_1A_2)$ and 
$\tr(A_1)
\tr(A_2)$ (corresponding to $c_2$ and $c_1^2$, respectively).  

For $M$, $\tr(A_1A_2)$ gives rise 
to the classical Chern-Simons invariant for $M$.  However, the Chern-Simons class associated to 
$\tr(A_1)\tr(A_2)$ involves $\tr(\omega_1-\omega_0)\tr(\Omega_t)$, 
which vanishes since both forms take values in skew-symmetric endomorphisms.

In contrast, on $LM$ we know that the WCS class $CS^W_3$ associated to 
$\tr(A_1A_2)$ vanishes.  The WCS associated to $\tr(A_1)\tr(A_2)$ involves 
$\tr\sigma_{-1}(\omega_1-\omega_0) = \tr\sigma_{-1}(\omega_1)$ and $\tr\sigma_{-1}(\Omega_t).$  
Both $\omega_1$ and $ \Omega_t$ take values in skew-symmetric $\pdo$s, but
this does not imply that the terms in their symbol expansions are skew-symmetric.  In fact, a calculation using Lemma \ref{old2.1} shows that $\sigma_{-1}(\omega_1)$ is not skew-symmetric. 
Thus the WCS class associated to $\tr(A_1)\tr(A_2)$ may be nonzero.

\end{exm}

\section{{\bf An Application of Wodzicki-Chern-Simons Classes to Circle
    Actions}}\label{dimfive}

In this section we use WCS classes to distinguish different $S^1$ actions on
$M=S^2\times S^3$. We use this to conclude that $\pi_1(\diffm, id)$ is infinite.  

Recall that  $H^*(LM)$ denotes de Rham cohomology of complex valued 
forms.   In particular, integration of closed forms over homology cycles gives a pairing of
$H^*(LM)$ and $H_*(LM,\BbC)$. 

 For any closed oriented manifold $M$, let $a_0,a_1:S^1\times M\to M$ be two smooth actions.  Thus
$$a_i(0,m) = m, \ a_i(\theta,a(\psi,m)) = a_i(\theta + \psi, m).$$

\begin{defn} (i)  $a_0$ and $a_1$
 are {\it smoothly  homotopic} if there exists a smooth map 
$$F:[0,1]\times S^1\times M\to M,\ F(0,\theta,m) = a_0(\theta,m),\
 F(1,\theta,m) = a_1(\theta,m).$$

(ii)  $a_0$ and $a_1$ are {\it smoothly  homotopic through actions} if
 $F(t,\cdot,\cdot):S^1\times M\to M$ is an action for all $t$.

\end{defn}

We can rewrite an action in two equivalent ways.

\begin{itemize}
\item $a$ determines (and is determined by)
$a^D:S^1\to \diff(M)$ given by
$a^D(\theta)(m) = a(\theta,m).$  $a^D(\theta)$ is a diffeomorphism because 
$$a^D(-\theta)(a^D(\theta,m)) = a(-\theta, a(\theta,m)) = m.$$
Since $a^D(0) = id,$  we get a class $[a^D]\in \pi_1(\diff(M), id)$, the
fundamental group of $\diff(M)$ based at $id.$  Here Diff$(M)$ is a Banach manifold
as an open subset of the Banach manifold of $\maps(M) = \maps(M,M)$ of some fixed Sobolev class.

\item $a$ determines (and is determined by)
$a^L:M\to LM$ given by $ a^L(m)(\theta) = a(\theta,m)$. This determines a class
$[a^L]\in H_n(LM,\Z)$ with $n = {\rm dim}(M)$ by setting $[a^L] = a^L_*[M].$
  In concrete terms, if we triangulate $M$ as the $n$-cycle $\sum_i n_i\sigma_i$,
with $\sigma_i:\Delta^n\to M$, 
then $[a^L]$ is the homology class of
the cycle $\sum_i n_i (a^L\circ \sigma_i).$  

\end{itemize}

We give a series of elementary lemmas comparing these maps.

\begin{lem}\label{lem:one} $a_0$ is smoothly homotopic to $a_1$ through actions iff $[a^D_0] =
	  [a^D_1]\in \pi_1(\diff(M), id).$
\end{lem}

\begin{proof} ($\Rightarrow$) Given $F$ as above, set $G:[0,1]\times S^1\to
	  \diff(M)$ by $G(t,\theta)(m) = F(t,\theta,m).$  We have $G(0,\theta)(m)
	  = a_0(\theta,m) = a^D(\theta)(m)$, $G(1,\theta)(m) = a_1(\theta, m)
	  = a^D_1(\theta)(m)$.
$G(t,\theta)\in\diff(M)$, because 
$$G(t,-\theta)(G(t,\theta)(m)) = F(t,-\theta,F(t,\theta,m)) = F(t,0,m) = m.$$
(This uses that $F(t,\cdot,\cdot)$ is an action.)
Since $F$ is
smooth, $G$ is a continuous (in fact, smooth) map of $\diff(M)$.
Thus $a^D_0, a^D_1$ 
are homotopic as elements of\\
 $\maps(S^1,\diff(M))$, so $[a^D_0] = [a^D_1].$
\bigskip

\noindent ($\Leftarrow$)  Let $G:[0,1]\times S^1\to \diff(M)$ be a continuous
homotopy from
$a^D_0(\theta) = G(0,\theta) $ to $a^D_1(\theta) = G(1,\theta)$ with $G(t,0) = id$
for all $t$. 
It is possible to approximate
$G$ arbitrarily well by a smooth map, since $[0,1]\times S^1$ is compact.  Set
$F:[0,1]\times S^1\times M\to M$ by
$F(t,\theta,m) = G(t,\theta)(m).$  
 $F$ is smooth.  Note that
$F(0,\theta,m) = G(t,\theta)(m) = a^D_0(\theta)(m) = a_0(\theta,m)$, and
$F(1,\theta,m) = a_1(\theta,m).$  Thus $a_0$ and $a_1$ are smoothly homotopic.
\end{proof}

%

There are similar results for $a^L.$

\begin{lem}\label{lem:three}  $a_0$ is smoothly homotopic to $a_1$ iff $a^L_0,
a^L_1:M\to LM$ are smoothly homotopic.
\end{lem}

\begin{proof}  Let $F$ be the homotopy from $a_0$ to $a_1$.  Set
	  $H:[0,1]\times M \to LM $ by $H(t,m)(\theta) = F(t,\theta,m).$  Then
$H(0,m)(\theta) = F(0,\theta,m) = a_0(\theta,m) =  a^L_0(m)(\theta)$,
  $H(1,m)(\theta) =  a^L_1(m)(\theta),$  so $H$ is a homotopy from $
  a^L_0$ to $ a^L_1.$  
It is  easy to check that $H$ is smooth.

Conversely, if $H:[0,1]\times M \to LM $ is a smooth homotopy from $a^L_0$ to
$a^L_1$, set $F(t,\theta, m) = H(t,m)(\theta).$  
\end{proof}

\begin{cor}\label{cor:one} If $a_0$ is smoothly homotopic to $a_1$, then
$[a^L_0] =  [a^L_1]\in H_n(LM,\Z).$  
\end{cor}

\begin{proof} By the last Lemma, $a^L_0$ and $a^L_1$ are homotopic. Thus 
$[a^L_0] = a^L_{0,*}[M] = a^L_{1,*}[M] = [a^L_1].$
\end{proof}

This yields a technique to use WCS classes to distinguish actions and to investigate 
$\pi_1(\diffm,id).$  From now on, ``homotopic" means ``smoothly homotopic."

\begin{prop} \label{prop:two} Let dim$(M)=2k-1.$ Let $a_0, a_1:S^1\times M\to M$ be actions.


(i)  If $\int_{[a^L_0]} CS^W_{\kk} \neq \int_{[a^L_1]} CS^W_{\kk}$, then $a_0$ and $a_1$
  are not homotopic through actions, and $[a^D_0]\neq [a^D_1]\in \pi_1(\diff(M),id).$

(ii)  If $\int_{[a_1^L]} CS^W_{\kk} \neq 0,$ then
  $\pi_1(\diff(M), id)$  is infinite.

\end{prop}

\begin{proof} 

(i) By Stokes' Theorem, $[a^L_0]\neq [a^L_1]\in H_n(LM,\C)$.
  By Corollary \ref{cor:one}, $a_0$ and $a_1$ are not homotopic,
   and hence not homotopic
  through actions.  By
Lemma \ref{lem:one}, $[a^D_0]\neq [a^D_1]\in \pi_1(\diff(M),id).$

(ii) Let $a_n$ be the $n^{\rm th}$ iterate of 
 $a_1$, i.e. $a_n(\theta,m) =
a_1(n\theta,m).$  

We claim that 
 $\int_{[a^L_n]}CS^W_{\kk} =
n\int_{[a^L_1]}CS^W_{\kk}$.  By (\ref{5.4}), every term in $CS^W_{\kk}$ is of the
form $\ints\dot\gamma(\theta) f(\theta)$, where $f$ is a periodic function on the
circle.  Each loop $\gamma\in
a^L_1(M)$ corresponds to the loop $\gamma(n\cdot)\in a^L_n(M).$  Therefore the term
$\ints\dot\gamma(\theta) f(\theta)$ is replaced by 
$$\ints \frac{d}{d\theta}\gamma(n\theta) f(n\theta)d\theta 
 = n\int_0^{2\pi} \dot\gamma(\theta)f(\theta)d\theta.$$
Thus $\int_{[a^L_n]}CS^W_{\kk} = n\int_{[a^L_1]}CS^W_{\kk}.$
 By (i), the $[a^L_n]\in 
\pi_1(\diff(M), id)$
are all distinct.  

\end{proof}

\begin{rem}
If two actions
are homotopic through actions,
the $S^1$ index of an equivariant operator of the two actions is the same. (Here equivariance
means for each action $a_t, t\in [0,1].$)
In contrast to Proposition \ref{prop:two}(ii), the $S^1$ index of an equivariant operator
cannot distinguish actions on odd dimensional manifolds, as the
$S^1$ index vanishes. This can be seen from the
local version of the
$S^1$ index theorem \cite[Thm. 6.16]{BGV}. For the normal bundle to the
fixed point set is always even dimensional, so the fixed point set consists of
odd dimensional submanifolds.  The integrand in the fixed point submanifold
contribution to the $S^1$-index is the constant term in the short time
asymptotics of the appropriate heat kernel.  In odd dimensions, this constant
term is zero.

  In \cite{MRT2}, we  interpret the $S^1$ index theorem as
the integral of an equivariant characteristic class over $[a^L]$.
\end{rem}

We now apply these methods to a Sasaki-Einstein metric on $S^2\times S^3$
constructed in \cite{gdsw}
 to prove the following:

\begin{thm} (i) There is an $S^1$ action on $S^2\times S^3$ that is not smoothly homotopic
to the trivial action.

(ii) $\pi_1(\diff(S^2\times S^3), id)$ is infinite.
\end{thm}

The content of  (i) is that although the $S^1$-orbit
 $\gamma_x$ through $x\in S^2\times S^3$
is contractible to $x$, the contraction cannot be constructed to be 
smooth in $x$.  

\begin{proof}
According to  \cite{gdsw}, the locally defined metric 
\begin{eqnarray}\label{metric} g &=&\frac{1-cy}{6}(d\theta^2 + \sin^2\theta d\phi^2) + 
\frac{1}{w(y)q(y)} dy^2 + \frac{q(y)}{9}[d\psi^2 -\cos\theta d\phi^2]\nonumber\\
&&\qquad + w(y)\left[d\alpha + \frac{ac-2y+y^2c}{6(a-y^2)}[d\psi -\cos\theta d\phi]\right]^2,
\end{eqnarray}
with 
$$w(y) = \frac{2(a-y^2)}{1-cy}, q(y) = \frac{a-3y^2+2cy^3}{a-y^2},$$
is a family of Sasaki-Einstein metrics on a coordinate ball in the variables
$(\phi, \theta, \psi, y, \alpha).$  Here $a$ and $c$ are constants, and we can take $a\in (0,1], c=1$.  
For $p,q$ relatively prime, $q<p$, and satisfying $4p^2-3q^2 = n^2$ for some integer $n$,  
and for  $a = a(p,q)< 1$, the metric extends to a $5$-manifold $Y^{p,q}$ which has
the coordinate ball as a dense subset.  
In this case, $( \phi, \theta, \psi, y)$
are  spherical coordinates on $S^2\times S^2$ with a nonstandard metric, and $\alpha$ is the fiber coordinate of an 
$S^1$-fibration $S^1\to Y^{p,q}\to S^2\times S^2.$  
$Y^{p,q}$ is diffeomorphic to $S^2\times S^3,$ and
 has first Chern class which integrates over the two $S^2$ factors to $p+q$ and $p$ \cite[\S2]{gdsw}.
The coordinate ranges are $\phi\in (0,2\pi), \theta \in (0,\pi), \psi\in (0,2\pi)$,
$\alpha\in (0,2\pi\ell)$,
where $\ell = \ell(p,q)$, and 
$y\in (y_1, y_2)$, with the $y_i$ the two smaller roots of $a-3y^2+2y^3=0$.  $p$ and
$q$ determine $a, \ell, y_1, y_2$ explicitly \cite[(3.1), (3.4), (3.5), (3.6)]{gdsw}.

For these choices of $p, q$, we get an $S^1$-action $a_1$ on $Y^{p,q}$ by rotation in the $\alpha$-fiber.
We claim that for e.g. $(p,q) = (7,3)$,
\begin{equation}\label{neqo} \int_{[a_1^L]} CS_5^W(g) \neq 0.
\end{equation}
By Proposition \ref{prop:two}(iii), this implies $\pi_1(\diff(S^2\times S^3), id)$ is infinite.
Since the trivial action $a_0$ has $\int_{[a_0^L]} CS_5^W(g) = 0$ (by the proof of 
Proposition \ref{prop:two}(ii) with $n=0$), $a_0$ and $a_1$ are not smoothly homotopic by
Proposition \ref{prop:two}(i).   Thus showing (\ref{neqo}) will prove the theorem.

Set $M = S^2\times S^3$. Since $a_1^L:M\to LM$ has degree one on its image,
\begin{equation}\label{cswf} 
\int_{[a_1^L]} CS_5^W(g) = \int_M a_1^{L,*} CS_5^W(g).
\end{equation}
For $m\in M$, 
$$a_1^{L,*}CS_5^W(g)_m = f(m)d\phi\wedge d\theta\wedge dy\wedge d\psi\wedge 
d\alpha$$ 
for some $f\in C^\infty(M)$.  We determine
$f(m)$ by explicitly computing $a_{1,*}^L(\partial_\phi),..., a_{1,*}^L(\partial_\alpha),$ (e.g. 
$a_{1,*}^L(\partial_\phi)(a^L(m))(t) = \partial_\phi|_{a(m,t)}$ ),
and noting
\begin{eqnarray}\label{f} f(m) &=& f(m)d\phi\wedge d\theta\wedge dy\wedge d\psi\wedge 
d\alpha(\partial_\phi,\partial_\theta,\partial_y,\partial_\psi,\partial_\alpha)\nonumber\\
&=& a_1^{L,*}CS_5^W(g)_m(\partial_\phi,...,\partial_\alpha)\\
&=& CS_5^W(g)_{a_1^L(m)}(a_{1,*}^L(\partial_\phi),...,a_{1,*}^L(\partial_\alpha)).\nonumber
\end{eqnarray}
Since $CS^W_5(g)$ is explicitly computable from the formulas in \S\ref{localsymbols}, we can
compute $f(m)$ from (\ref{f}).  Then $\int_{[a_1^L]} CS_5^W(g) = \int_{M} f(m)
d\phi\wedge d\theta\wedge dy\wedge d\psi\wedge 
d\alpha$ can be computed as an ordinary integral in the dense coordinate space.

Via this method, in the Mathematica file {\tt ComputationsChernSimonsS2xS3.pdf} at {\tt http://math.bu.edu/people/sr/},  $\int_{[a_1^L]} CS_5^W(g)$ is computed as
 a function of $(p,q).$  For example, 
$(p,q) = (7,3)$,
$$\int_{[a_1^L]} CS_5^W(g)  = -\frac{1849\pi^4}{22050}.$$
This formula is exact; the rationality up to $\pi^4$ follows from $4p^2-3q^2$ being a perfect
square, as then the various integrals computed in (\ref{cswf}) with respect to our coordinates
are rational functions evaluated at rational endpoints. 
 In particular, (\ref{neqo})
holds.
  \end{proof}

\begin{rem}  
For $a=1$, the metric extends to the closure of the coordinate chart, but the total space is $S^5$ with the standard metric.  
$\pi_1(\diff(S^5))$ is torsion \cite{F-H}.  By Proposition \ref{prop:two}(ii), $\int_{[a^L]} CS^W_5 = 0$
for any circle action on $S^5.$  In the formulas in the Mathematica file, $\int_{[a^L]}CS^W_5$ is proportional to
$(-1+a)^2$, which vanishes at $a=1$.  
This gives a check of the validity of the computation.
\end{rem}

\bibliographystyle{amsplain}
\bibliography{Paper}

\bigskip
\hfill \today \\

\end{document}